\pdfoutput=1
\RequirePackage{ifpdf}
\ifpdf
\documentclass[pdftex]{sigma}
\else
\documentclass{sigma}
\fi

\numberwithin{equation}{section}

\newtheorem{thm}{Theorem}[section]
\newtheorem{cor}[thm]{Corollary}
\newtheorem{lem}[thm]{Lemma}
\newtheorem{prop}[thm]{Proposition}
{\theoremstyle{definition}
\newtheorem{exam}[thm]{Example}
\newtheorem{defn}[thm]{Definition}
}

\begin{document}

\allowdisplaybreaks

\renewcommand{\PaperNumber}{019}

\FirstPageHeading

\ShortArticleName{Tight Frame with Hahn and Krawtchouk Polynomials of Several Variables}

\ArticleName{Tight Frame with Hahn and Krawtchouk Polynomials\\
of Several Variables}

\Author{Yuan XU}

\AuthorNameForHeading{Y.~Xu}

\Address{Department of Mathematics, University of Oregon, Eugene, Oregon 97403-1222, USA}
\Email{\href{mailto:yuan@uoregon.edu}{yuan@uoregon.edu}}
\URLaddress{\url{http://pages.uoregon.edu/yuan/}}

\ArticleDates{Received November 06, 2013, in f\/inal form February 25, 2014; Published online March 03, 2014}

\Abstract{Finite tight frames for polynomial subspaces are constructed using monic Hahn polynomials and Krawtchouk
polynomials of several variables.
Based on these polynomial frames, two methods for constructing tight frames for the Euclidean spaces are designed.
With ${\mathsf r}(d,n):= \binom{n+d-1}{n}$, the f\/irst method generates, for each $m \ge n$,
two families of tight frames in ${\mathbb R}^{{\mathsf r}(d,n)}$ with ${\mathsf r}(d+1,m)$ elements.
The second method generates a~tight frame in~${\mathbb R}^{{\mathsf r}(d,N)}$ with $1 + N \times{\mathsf r}(d+1, N)$ vectors.
All frame elements are given in explicit formulas.}


\Classification{33C50; 42C15}

\section{Introduction}

For a~real f\/inite dimension space $V$ with inner product $(\cdot,\cdot)$, a~set of vectors $\{v_k: k =1,\ldots, M\}$
in $V$ is called a~(Parseval) tight frame with $M$ elements if
\begin{gather*}
(v, v) =\sum\limits_{k=1}^M (v, v_k)^2,
\qquad
\forall \, v \in V.
\end{gather*}
Tight frames are connected to a~wide variety of problems in pure and applied mathematics and they have been studied
intensively in recent years (see, for examples,~\cite{BE, BF, Cass, CK, Heil} and references therein).
We construct tight frames with the Hahn and Krawtchouk polynomials of several variables, and use them to construct
several families of tight frames for the Euclidean spaces that are given in explicit formulas.

Let $\Pi_n^d$ denote the space of polynomials of total degree at most $n$ in $d$-variables.
Let ${\langle} \cdot, \cdot{\rangle}$ denote an inner product def\/ined on the space of polynomials or on $\Pi_N^d$ for
a~f\/ixed integer $N$.
A polynomial $p \in \Pi_n^d$ is called an orthogonal polynomial of degree $n$ with respect to the inner product if
\begin{gather*}
{\langle} p, q {\rangle} = 0
\qquad
\text{for all}
\quad
q \in \Pi_{n-1}^d.
\end{gather*}
Let ${\mathcal V}_n^d$ denote the space of orthogonal polynomials of degree $n$.
Then
\begin{gather*}
\operatorname{dim}{\mathcal V}_n^d = \binom{n+d-1}{n},
\qquad
n = 0,1,2, \ldots.
\end{gather*}
For $\nu \in {\mathbb N}_0^d$, let $|\nu|:= \nu_1+\dots + \nu_d$.
A basis $\{P_\nu: |\nu| = n \}$ of ${\mathcal V}_n^d$ is called mutually orthogonal if ${\langle} P_\nu, P_\mu {\rangle}
= 0$ whenever $\nu \ne \mu$ and $|\nu| = |\mu| =n$, and it is called orthonormal if, in addition, ${\langle} P_\nu,
P_\nu {\rangle} = 1$.
The space ${\mathcal V}_n^d$ can have many distinct bases.

For $a \in {\mathbb R}$ and $n \in {\mathbb N}_0^d$, let $(a)_n = a(a+1)\cdots (a+n-1)$ be the Pochhammer symbol.
For $x \in {\mathbb R}^d$ and $\nu \in {\mathbb N}_0^d$ we let $x^\nu: = x_1^{\nu_1} \cdots x_d^{\nu_d}$, and let
\begin{gather*}
\nu !:= \nu_1! \cdots \nu_d!
\qquad
\text{and}
\qquad
(\nu)_\mu:= (\nu_1)_{\mu_1} \cdots (\nu_d)_{\mu_d}.
\end{gather*}
The Hahn polynomials and Krawtchouk polynomials of $d$-variables are discrete orthogonal polynomials def\/ined on the
homogeneous coordinates of
\begin{gather*}
{\mathbb Z}_N^{d+1}:= \big\{ {\alpha} \in {\mathbb Z}^{d+1}: |{\alpha}| = {\alpha}_1 + \dots + {\alpha}_{d+1} = N \big\}.
\end{gather*}
For ${\kappa} \in {\mathbb R}^{d+1}$ with ${\kappa}_1 > -1, \ldots, {\kappa}_{d+1} > -1$, the Hahn polynomials are
orthogonal with respect to the inner product
\begin{gather*}
{\langle} f, g{\rangle}_{{\mathsf H}_{{\kappa},N}}: = \frac{N!}{(|{\kappa}|+d+1)_N}
\sum\limits_{x \in Z_N^{d+1}} f(x) g(x)
{\mathsf H}_{\kappa}(x),
\qquad
{\mathsf H}_{\kappa}(x): = \frac{({\kappa}+{\mathbf{1}})_{\alpha}}{{\alpha}!},
\end{gather*}
where ${\mathbf{1}}: = (1,\ldots, 1)$.
For $\rho \in {\mathbb R}^{d}$ with $0 < \rho_1 , \ldots, \rho_d < 1$ and $|\rho| <1$, the Krawtchouk polynomials are
orthogonal with respect to the inner product
\begin{gather*}
{\langle} f, g{\rangle}_{{\mathsf K}_{\rho,N}}: = \frac1{N!}\sum\limits_{x \in Z_N^{d+1}} f(x) g(x){\mathsf K}_\rho(x),
\qquad
{\mathsf K}_\rho(x): = \frac{{\boldsymbol{\rho}}^{\kappa}}{x!},
\end{gather*}
where ${\boldsymbol{\rho}} = (\rho, 1-|\rho|)$.
These polynomials have been studied extensively in the literature.
For this paper, the most relevant work on the Hahn polynomials are those in~\cite{GS11,GS13, IX07, KM, HR, Tr1, Tr2, Tr3, X13}.
For the Krawtchouk polynomials, we refer to~\cite{DG, Tr2, X13} and references therein.
In both cases, a~family of mutually orthogonal polynomials can be def\/ined in terms of the classical Hahn polynomials or
Krawtchouk polynomials of one variable, and they have been used in most of the studies in the literature.

In the present paper, we study monic Hahn and Krawtchouk polynomials instead.
These are orthogonal projections of $(-1)^{|{\alpha}|}(-x)_{\alpha}$ on either ${\mathcal V}_n^d({\mathsf H}_{{\kappa},N}
)$ or ${\mathcal V}_n^d({\mathsf K}_{\rho,N}
)$.
For ${\alpha} \in {\mathbb Z}_m^{d+1}$ and $m \ge n$, let ${\mathsf U}_{{\alpha},n}
(x) =(-1)^{|{\alpha}|}(-x)_{\alpha} + \cdots$ denote the monic orthogonal polynomial of degree~$n$.
Then $\Xi_{n,m}:= \{
{\mathsf U}_{{\alpha},n}
: {\alpha} \in {\mathbb Z}_m^{d+1}\}$ is a~subset of ${\mathcal V}_n^d$ and it contains an orthogonal basis of~${\mathcal V}_n^d$.
In the special case of $m=n$ and ${\alpha} = ({\alpha}', n-|{\alpha}'|)$ with ${\alpha}' \in {\mathbb N}_0^d$ and
$|{\alpha}'| \le n$, the monic Hahn and Krawtchouk polynomials were studied in~\cite{Tr1, Tr2}.
We shall describe generating functions, connection to the mutually orthogonal polynomials, and explicit expressions for
these polynomials.
Furthermore, the set $\Xi_{n,m}$ contains far more elements than the dimension of ${\mathcal V}_n^d$, we shall prove
that the set $\Xi_{n,m}$ is a~tight frame for the space ${\mathcal V}_n^d$ and these frames are given in explicit
formulas.
The connection between tight frames and orthogonal polynomials of several variables has been explored for the Jacobi
polynomials on the simplex \cite{W} (see also~\cite{W2, X13}).
The discrete orthogonal polynomials have not been used in this connection until now.

A f\/inite tight frame is equivalent to a~tight frame on an Euclidean space.
From our Hahn and Krawtchouk polynomial frames, we are able to derive a~large family of tight frames for Euclidean
spaces that are given by explicit formulas.
There are many ways of constructing tight frames for the Euclidean spaces.
One can, for example, use the classical discrete Fourier transform to construct tight frame for ${\mathbb C}^n$ or
${\mathbb R}^n$ that contains $m$ elements for all $m > n$.
There are also algorithms that can be used to construct tight frames with certain additional features~\cite{Cass}.
Our frames are all given by explicit formulas.
In fact, we shall give two methods for constructing such frames.
For $n, d \in {\mathbb N}$, let
\begin{gather*}
{\mathsf r} (d, n):= \binom{n+d-1}{n} = \operatorname{dim}{\mathcal V}_n^d
\qquad
\text{and}
\qquad
{\mathsf n}(d,n): = {\mathsf r}(d+1,n).
\end{gather*}
For a~positive integer $m \ge n$, our f\/irst method uses either the Hahn polynomial frames or the Krawtchouk polynomial
frames to construct explicitly tight frames in ${\mathbb R}^{{\mathsf n}(d,n)}$ that contain ${\mathsf n}(d,m)$ elements.
Our second method uses the Hahn polynomial frames with the parameter ${\kappa}=0$ to construct tight frames in
${\mathbb R}^{{\mathsf n}(d,N)}$ with $N \times{\mathsf n}(d,N)+1$ elements.
All frames are given by explicit formulas and those derived from the second method have the additional feature that all
but one frame elements are vectors whose entries consist of only rational numbers.

The paper is organized as follows.
In the next section, we recall necessary background on orthogonal polynomials.
Monic Hahn polynomials are introduced and developed in Section~\ref{Sec3}, from which monic Krawtchouk polynomials are deduced by
taking an appropriate limit in Section~\ref{Sec4}.
Tight frames for the Hahn and the Krawtchouk polynomials, as well as for the Euclidean spaces are discussed
in Section~\ref{Sec5}, where several examples are explicitly given.

\section{Background and preliminary}

Throughout this paper we reserve the Greek letter ${\kappa}$ and $\rho$ for the parameters in the weight function ${\mathsf H}_{{\kappa},N}$
and ${\mathsf K}_\rho$.
We shall use other Greek letters for the multiindexes.
Orthogonal polynomials in ${\mathcal V}_n^d$ are indexed by $\nu$ or $\mu$ in ${\mathbb N}_0^d$, and we use
${\alpha}$, ${\beta}$, ${\gamma}$ for indexes in ${\mathbb N}_0^{d+1}$.

For $a, b > -1$, the classical Hahn polynomial $Q_n(x; a,b,N)$ is a~${}_3F_2$ hypergeometric function given by
\begin{gather}
\label{eq:HahnQ}
{\mathsf Q}_n(x; a, b, N):= {}
_3 F_2 \left(
\begin{matrix}
-n, n+a + b+1, -x
\\
a+1, -N
\end{matrix}
; 1\right),
\qquad
n = 0, 1, \ldots, N.
\end{gather}
These are discrete orthogonal polynomials def\/ined on the set $\{0,1,\ldots, N\}$.
For the Hahn polynomials of $d$-variables, a~mutually orthogonal basis for ${\mathcal V}_n^d({\mathsf H}_{{\kappa},N}
)$ can be given in terms of ${\mathsf Q}_n(x; a,b,N)$.
We need the following def\/inition.

For $y=(y_1,\ldots, y_{d}) \in {\mathbb R}^{d}$ and $1 \le j \le d$, we def\/ine
\begin{gather}
\label{xsupj}
{\mathbf y}_j:= (y_1, \ldots, y_j)
\qquad
\text{and}
\qquad
{\mathbf y}^j:= (y_j, \ldots, y_d).
\end{gather}
We also def\/ine ${\mathbf y}_0:= \varnothing$ and ${\mathbf y}^{d+1}:= \varnothing$.
In particular, ${\mathbf y}_d = {\mathbf y}^1 = y$ and
\begin{gather*}
|{\mathbf y}_j| = y_1 + \dots + y_j,
\qquad
|{\mathbf y}^j| = y_j + \dots + y_d,
\qquad
\text{and}
\qquad
|{\mathbf y}_0| = |{\mathbf y}^{d+1}| = 0.
\end{gather*}
For the parameter vector ${\kappa}=({\kappa}_1,\ldots, {\kappa}_{d+1})$, we def\/ined ${\kappa}^j:= ({\kappa}_j, \ldots,
{\kappa}_{d+1})$ for $1 \le j \le d+1$.
For $\nu \in {\mathbb N}_0^d$ and $x \in {\mathbb R}^d$, def\/ine
\begin{gather}
\label{eq:aj}
a_j:=a_j(\kappa,\nu):=|\kappa^{j+1}| + 2 |\nu^{j+1}| + d-j,
\qquad
1 \le j \le d.
\end{gather}

\begin{prop}
For $x \in {\mathbb Z}_N^{d+1}$ and $\nu \in {\mathbb N}_0^d$, $|\nu| \le N$, define
\begin{gather}
{\mathsf H}_\nu(x;\kappa, N):=\frac{(-1)^{|\nu|}}{(-N)_{|\nu|}}
\prod\limits_{j=1}^d \frac{(\kappa_j+1)_{\nu_j}}{(a_j+1)_{\nu_j}}\big(-N+|{\mathbf x}_{j-1}|+|\nu^{j+1}|\big)_{\nu_j}
\nonumber
\\
\phantom{{\mathsf H}_\nu(x;\kappa, N):=}
\times
{\mathsf Q}_{\nu_j}\big(x_j; \kappa_j, a_j, -N+|{\mathbf x}_{j-1}|-|\nu^{j+1}|\big).
\label{eq:Hn-prod}
\end{gather}
Then the polynomials in $\{{\mathsf H}_\nu(x; \kappa,N): |\nu| = n\}$
form a~mutually orthogonal basis of ${\mathcal V}_n^d({\mathsf H}_{{\kappa},N})$
and ${\mathsf B}_\nu:= {\langle}{\mathsf H}_\nu(\cdot; \kappa,N), \;{\mathsf H}_\nu(\cdot; \kappa,N) {\rangle}_{{\mathsf H}_{{\kappa},N}}$
is given by, setting ${\lambda}_{\kappa}: = |{\kappa}|+d+1$,
\begin{gather}
\label{eq:Bnu}
{\mathsf B}_\nu({\kappa}, N):=\frac{(-1)^{|\nu|}
({\lambda}_k)_{N+|\nu|}}{(-N)_{|\nu|} ({\lambda}_k)_N ({\lambda}_k)_{2|\nu|}} \prod\limits_{j=1}^d
\frac{({\kappa}_j+a_j+1)_{2\nu_j} ({\kappa}_j+1)_{\nu_j} \nu_j! }{ (\kappa_j+a_j + 1)_{\nu_j} (a_j+1)_{\nu_j}}.
\end{gather}
\end{prop}

This basis was def\/ined in~\cite{KM} and given in the present form in~\cite{IX07,Tr2, X05}.
Another basis for ${\mathcal V}_n^d({\mathsf H}_{{\kappa},N})$
that is of interests for our work is the monic orthogonal basis studied in~\cite{Tr1}, which will come out as
a~special case of the discussion in the next section.
It is worth mentioning that, for $\nu \in {\mathbb N}_0^d$, $m,N = 0,1,\dots$,
\begin{gather}
\label{eq:B-A}
\frac{
{\mathsf B}_\nu({\kappa},N)}{
{\mathsf B}_\nu({\kappa},m)}
= \frac{(|{\kappa}|+d+1)_{N+|\nu|}(|{\kappa}|+d+1)_{m} (-m)_{|\nu|}}{(|{\kappa}|+d+1)_{m+|\nu|}(|{\kappa}|+d+1)_{N}(-N)_{|\nu|}}
\end{gather}
depends only on $|\nu|$, not the elements of $\nu$.
The expansion of ${\mathsf H}_\nu(x;{\kappa},N)$ in the shifted monic monomials $(-x)_{\alpha}$ is given in~\cite{X13}.

\begin{prop}
For $\nu \in {\mathbb N}_0^d$ with $|\nu| =n$ and $x \in {\mathbb Z}_N^{d+1}$,
\begin{gather}
\label{eq:phi-monic}
{\mathsf H}_\nu(x; {\kappa}, N) = \frac{n!}{(-N)_n}
\sum\limits_{|{\alpha}|=n}
\frac{
{\mathsf H}_\nu({\alpha}; {\kappa}, n)}{{\alpha}!} (-x)_{\alpha}.
\end{gather}
\end{prop}

There is a~close relation between the Hahn polynomials and the orthogonal polynomials with respect to the weight
function
\begin{gather*}
W_{\kappa}(x):= x_1^{{\kappa}_1}\cdots x_d^{{\kappa}_d} (1-|x|)^{{\kappa}_{d+1}},
\qquad
{\kappa}_1 > -1,\ldots, {\kappa}_{d+1} > -1,
\end{gather*}
on the simplex $T^d:= \{x\in {\mathbb R}^d: x_1 \ge 0, \ldots, x_d \ge 0, |x| \le 1\}$.
Let
\begin{gather*}
{\langle} f, g {\rangle}_{W_{\kappa}}:= w_{\kappa} \int_{T^d} f(x) g(x) W_{\kappa}(x) dx
\qquad
\text{with}
\qquad
w_\kappa:= \frac{\Gamma(|\kappa| + d+1)}{\prod\limits_{i=1}^{d+1}\Gamma(\kappa_i +1)}.   
\end{gather*}
A mutually orthogonal basis for ${\mathcal V}_n^d(W_{\kappa})$ can be given in terms of the classical Jacboi polynomials
\begin{gather*}
\frac{P_n^{(a,b)}(t)}{P_n^{(a,b)}(1)} = {}_2 F_1\left(-n,n+a+b+1; a+1; \frac{1-t}{2}\right).
\end{gather*}

\begin{prop}
For $\nu \in {\mathbb N}_0^d$ and $x \in {\mathbb R}^d$, define
\begin{gather*}
P_\nu^\kappa (x):= \prod\limits_{j=1}^d \left(\frac{1-|{\mathbf x}_j|}{1-|{\mathbf x}_{j-1}|} \right)^{|\nu^{j+1}|}
\frac{P_{\nu_j}^{(a_j,\kappa_j)}\left (\frac{2x_j}{1-|{\mathbf x}_{j-1}|} -1\right)}{P_{\nu_j}^{(a_j,\kappa_j)}(1)},
\end{gather*}
where $a_j$ is defined by~\eqref{eq:aj}.
Then the polynomials in $\{P_\nu^{\kappa}:|\nu|=n\}$ form a~mutually orthogonal basis of ${\mathcal V}_n^d(W_{\kappa})$
with $\langle P_\nu^{\kappa}, P_\nu^{\kappa} \rangle_{W_\kappa}$ given by
\begin{gather*}
\langle P_\nu^{\kappa}, P_\nu^{\kappa} \rangle_{W_\kappa} = \frac1{(|\kappa|+d+1)_{2|\nu|}} \prod\limits_{j=1}^d
\frac{({\kappa}_j+a_j+1)_{2 \nu_j} ({\kappa}_j +1)_{\nu_j} \nu_j!}{(\kappa_j+a_j+1)_{\nu_j}(a_j+1)_{\nu_j}}.
\end{gather*}
\end{prop}

This basis is well studied, but our normalization is dif\/ferent from the usual def\/inition (see~\cite[p.~47]{DX})
by $P_{\nu_j}^{(a_j,\kappa_j)}(1)$ in the denominator.
As def\/ined in~\cite{KM} and later recognized in~\cite{X05}, the polynomials $P_\nu^\kappa$ serve as generating function
for the Hahn polynomials.

\begin{prop}
Let ${\kappa}\in{\mathbb R}^{d+1}$ with $\kappa_i>-1$ and $N\in{\mathbb N}$.
For $\nu\in{\mathbb N}_0^d$, $|\nu|\le N$, and $y = (y', y_{d+1}) \in {\mathbb R}^{d+1}$,
\begin{gather}
\label{Hahngenfunc}
|y|^N P_\nu^{\kappa} \Big (\frac {y'}{|y|} \Big) =
\sum\limits_{\alpha \in {\mathbb Z}_N^{d+1}
} \frac{N!}{\alpha!}{\mathsf H}_\nu(\alpha; \kappa,N)y^\alpha.
\end{gather}
\end{prop}

The Hahn polynomials also appear as connecting coef\/f\/icients between $P_\nu^{\kappa}$ and monic ortho\-go\-nal polynomials
$R_{\alpha}^{\kappa}$ def\/ined by, for ${\alpha} \in {\mathbb N}_0^{d+1}$ and $x \in {\mathbb R}^d$,
\begin{gather*}
R_{\alpha}^{\kappa}(x):= X^{\alpha} + q_{\alpha}(x),
\qquad
q \in \Pi_{|{\alpha}|-1}^d,
\qquad
\text{where}
\qquad
X = (x,1-|x|).
\end{gather*}
The explicit formula of $R_{\alpha}(x)$ was derived in~\cite{X05}.

\begin{thm}
For ${\alpha} \in {\mathbb N}_0^{d+1}$ and $x \in {\mathbb R}^d$,
\begin{gather*}
R_{\alpha}^{\kappa}(x) =\frac{ (-1)^{|{\alpha}|} ({\kappa}+\mathbf{1})_{\alpha}}{(|{\kappa}|+d+|{\alpha}|)_{|{\alpha}|}}
\sum\limits_{{\gamma} \le {\alpha}}
\frac{(-{\alpha})_{{\gamma}} (|{\kappa}|+ d+|{\alpha}|)_{|{\gamma}|}}{({\kappa}+\mathbf{1})_{{\gamma}}{}{\gamma}!}X^{{\gamma}}.
\end{gather*}
Furthermore, $\{R_{\alpha}^{\kappa}: |{\alpha}| = n, {\alpha}_{d+1} = 0\}$ is a~basis of ${\mathcal V}_n^d(W_{\kappa})$.
\end{thm}

The cardinality of the set $\{R_{\alpha}^{\kappa}: {\alpha} \in {\mathbb N}_0^d, |{\alpha}| =N\}$ is larger than the
dimension of ${\mathcal V}_n^d(W_\kappa)$.
The Hahn polynomials serve as the connecting coef\/f\/icients of $P_\nu^{\kappa}$ and $R_{\alpha}^k$ as shown recently
in~\cite{X13}.

\begin{thm}
\label{thm:Pnu-Ra}
For $\nu \in {\mathbb N}_0^d$ with $|\nu|= n$ and ${\alpha}\in {\mathbb Z}_n^{d+1}$,
\begin{gather*}
P_\nu^{\kappa} (x) =\sum\limits_{|{\alpha}| =n}\frac{n!}{{\alpha}!}{\mathsf H}_\nu({\alpha}; {\kappa},n) R_{\alpha}^{\kappa}(x),
\end{gather*}
and, conversely,
\begin{gather}
\label{eq:Ra=Pnu}
R_{\alpha}^{\kappa}(x) = \frac{({\kappa}+{\mathbf{1}})_{\alpha}}{(|{\kappa}|+d+1)_n}
\sum\limits_{|\nu| =n}
\frac{
{\mathsf H}_\nu({\alpha}; {\kappa},n)}{
{\mathsf B}_\nu({\kappa},n)}
P_\nu^{\kappa}(x).
\end{gather}
\end{thm}

We now turn to basic results on the Krawtchouk polynomials.
For $0 < p < 1$, the classical Krawtchouk polynomial ${\mathsf K}_n(x; p, N)$ of one variable is def\/ined by
\begin{gather*}
{\mathsf K}_n(x; p, N):= {}_2 F_1 \left(
\begin{matrix}
-n, -x
\\
-N
\end{matrix}; \frac1{p}\right),
\qquad
n = 0, 1, \ldots, N.
\end{gather*}
For $\nu \in {\mathbb N}_0^d$, a~family of the Krawtchouk polynomials of $d$ variables can be given in terms of
Krawtchouk polynomials of one variable~\cite{IX07, Tr2}.
For $\rho \in {\mathbb N}_0^d$, we use the notation~\eqref{xsupj}, which implies that $|{\boldsymbol{\rho}}_j| =
\rho_1+\dots + \rho_j$ for $j =1,2,\ldots, d$.

\begin{prop}
Let $\rho \in {\mathbb R}^d$ with $0 < \rho_i < 1$ and $|\rho| < 1$.
For $\nu\in {\mathbb N}_0^d$, $|\nu|\le N$, and $x \in {\mathbb R}^d$, define
\begin{gather}
{\mathsf K}_\nu(x; \rho, N):=\frac{(-1)^{|\nu|}}{(-N)_{|\nu|}}
\prod\limits_{j=1}^d \frac{\rho_j^{\nu_j}}{(1- |{\boldsymbol{\rho}}_j|)^{\nu_j}} (-N+|{\mathbf x}_{j-1}|+|\nu^{j+1}|)_{\nu_j}
\nonumber\\
\phantom{{\mathsf K}_\nu(x; \rho, N):=}
\times
{\mathsf K}_{\nu_j}
\left(x_j; \frac{\rho_j}{1-|{\boldsymbol{\rho}}_{j-1}|}, N-|{\mathbf x}_{j-1}|- |\nu^{j+1}|\right).\label{eq:KrawK}
\end{gather}
The polynomials in $\{{\mathsf K}_\nu(\cdot;\rho, N): |\nu| = n\}$ form a~mutually
orthogonal basis of ${\mathcal V}_n^d({\mathsf K}_{\rho,N})$
and ${\mathsf C}_\nu (\rho, N): = {\langle}{\mathsf K}_\nu(\cdot;\rho, N), {\mathsf K}_\nu(\cdot;\rho, N){\rangle}_{{\mathsf K}_{\rho, N}}$
is given by
\begin{gather}
\label{eq:C(rho-N)}
{\mathsf C}_\nu (\rho, N): = \frac{(-1)^{|\nu|}
}{(-N)_{|\nu|}N!} \prod\limits_{j=1}^d \frac{\nu_j! \rho_j^{\nu_j} }{(1-|{\boldsymbol{\rho}}_j|)^{\nu_j-\nu_{j+1}} }.
\end{gather}
\end{prop}

The Krawtchouk polynomials in~\eqref{eq:KrawK} are limits of the Hahn polynomials in~\eqref{eq:Hn-prod}.
More precisely, setting ${\kappa} = t (\rho, 1-|\rho|)$, we have \cite{IX07}
\begin{gather}
\label{Hahn-Kraw}
\lim_{t \to \infty}{\mathsf H}_\nu(x; t(\rho, 1-|\rho|), N) = {\mathsf K}_\nu(x; \rho, N).
\end{gather}

\section{Monic Hahn polynomials of several variables}\label{Sec3}
\setcounter{equation}{0}

The monic polynomial $R_{\alpha}^{\kappa}$ has a~single monomial $x^{\alpha}$ as its highest term.
For the Hahn polynomials, the role of $x^{\alpha}$ is played by $\mathrm{m}_{\alpha}$ def\/ined as follows.

\begin{defn}
For $x \in {\mathbb R}^{d+1}$ and ${\alpha} \in {\mathbb N}_0^{d+1}$, define
\begin{gather*}
\mathrm{m}_{\alpha}(x):= (-1)^{|{\alpha}|} (-x)_{\alpha} = x^{\alpha} +q_{\alpha},
\qquad
q_{\alpha} \in \Pi_{n-1}^{d+1}.
\end{gather*}
\end{defn}

It follows that $\mathrm{m}_{\alpha}(x)$ is a~monic polynomial of degree $|{\alpha}|$.
For each ${\alpha} \in {\mathbb Z}_N^{d+1}$, we def\/ine the monic Hahn polynomial as the orthogonal polynomial that has
$\mathrm{m}_{\alpha}$ as its leading term.

\begin{defn}
For ${\alpha} \in {\mathbb N}_0^{d+1}$ and $|{\alpha}| \le N$, the monic Hahn polynomial
${\mathsf Q}_{\alpha}(\cdot; {\kappa}, N)$ in ${\mathcal V}_{|{\alpha}|}^d({\mathsf H}_{{\kappa},N})$ is defined uniquely by
\begin{gather*}
{\mathsf Q}_{\alpha}(x; {\kappa}, N):= \mathrm{m}_{\alpha}(x) + q_{\alpha}(x),
\qquad
q \in \Pi_{|{\alpha}|-1}^d,
\qquad
x \in {\mathbb Z}_N^{d+1}.
\end{gather*}
\end{defn}

The polynomials ${\mathsf Q}_{\alpha}(\cdot; {\kappa}, N)$ are discrete counterparts of monic orthogonal polynomial
$R_{\alpha}^{\kappa}$.
Like $R_{\alpha}$, the cardinality of the set $\{{\mathsf Q}_{\alpha}(\cdot; {\kappa}, N): |{\alpha}| =n\}\in {\mathcal V}_n^d(H_{{\kappa},N})$
is much larger than the dimension of ${\mathcal V}_n^d({\mathsf H}_{{\kappa},N})$,
so that the set contains redundancy.

In parallel to the relations between the monic orthogonal polynomials $R_{\alpha}^{\kappa}$ and the Jacobi polynomials
$P_\nu^{\kappa}$, in Theorem~\ref{thm:Pnu-Ra}, we can derive relations between
${\mathsf Q}_{\beta}(\cdot; {\kappa}, N)$ and the Hahn polynomials ${\mathsf H}_\nu(\cdot; {\kappa}, N)$.

\begin{prop}
\label{prop:H=Q=H}
For $x\in {\mathbb Z}_N^{d+1}$, ${\alpha} \in {\mathbb Z}_n^{d+1}$ and $\nu \in {\mathbb N}_0^d$ with $|\nu| =n$,
\begin{gather}
\label{eq:phi=xi}
{\mathsf H}_\nu(x;{\kappa},N) = \frac{(-1)^{n}
n!}{(-N)_n}
\sum\limits_{|{\alpha}| =n}
\frac{
{\mathsf H}_\nu({\alpha};{\kappa},n)}{{\alpha}!}{\mathsf Q}_{\alpha}(x;{\kappa},N),
\end{gather}
and, conversely,
\begin{gather}
\label{eq:xi=phi}
{\mathsf Q}_{\alpha}(x;{\kappa},N) =\frac{ ({\kappa}+{\mathbf{1}})_{\alpha} (-N)_n}{(-1)^n (|{\kappa}|+d+1)_n}
\sum\limits_{|\nu| =n}
\frac{
{\mathsf H}_\nu({\alpha};{\kappa},n)}{
{\mathsf B}_\nu({\kappa}, n)}{\mathsf H}_\nu(x;{\kappa},N).
\end{gather}
\end{prop}

\begin{proof}
Substituting $\mathrm{m}_{\beta}(x) =
{\mathsf Q}_{\beta} (x) - q_{\beta}$ into~\eqref{eq:phi-monic},
we see that~\eqref{eq:phi=xi} follows from the ortho\-go\-nality of ${\mathsf H}_\nu(\cdot; {\kappa}, N)$.
Conversely, since $\{{\mathsf H}_\nu(\cdot; {\kappa},N): |\nu| = n\}$ is a~basis of ${\mathcal V}_n^d({\mathsf H}_{{\kappa},N})$,
there are unique constants $c_{{\beta},\nu}$
such that ${\mathsf H}_\beta(\cdot; {\kappa}, N) = \sum\limits_{|\nu| = n}c_{{\beta},\nu}{\mathsf Q}_\nu(\cdot; {\kappa}, N)$.
Hence, by~\eqref{eq:phi=xi}, we obtain
\begin{gather*}
{\mathsf H}_\nu(x;{\kappa},N) = \frac{(-1)^{n}
n!}{(-N)_n}
\sum\limits_{|\mu| = n}
\sum\limits_{|{\beta}| =n}
\frac{
{\mathsf H}_\nu({\beta}; {\kappa},n)}{{\beta}!} c_{{\beta},\mu}{\mathsf H}_\mu(x;{\kappa},N).
\end{gather*}
Since ${\mathsf H}_\nu(\cdot, {\kappa}, N)$ are mutually orthogonal, we must have
\begin{gather*}
\frac{(-1)^{n}
n!}{(-N)_n}
\sum\limits_{|{\beta}| =n}
\frac{
{\mathsf H}_\nu({\beta}; {\kappa},n)}{{\beta}!} c_{{\beta},\mu} = \delta_{\nu,\mu}.
\end{gather*}
The orthogonality of ${\mathsf H}_\nu(\cdot;{\kappa},n)$ gives one solution of $c_{{\beta},\mu}$.
The uniqueness of $c_{{\beta},\mu}$ shows that it is the only solution and proves~\eqref{eq:xi=phi}.
\end{proof}

{\sloppy Throughout the rest of this section, we often use the abbreviation
\begin{gather*}
{\lambda}_k: = |{\kappa}| +d+1.
\end{gather*}
Recall that the Jacobi polynomials $P_\nu^{\kappa}$ are generating functions of the Hahn polynomials ${\mathsf H}_\nu(\cdot; {\kappa}, N)$.
It turns out that $R_{\alpha}^{\kappa}$ are generating functions of the monic orthogonal polynomials ${\mathsf
Q}_\nu(\cdot; {\kappa}, N)$.

}

\begin{thm}
Let $y = (y',y_{d+1}) \in {\mathbb R}^{d+1}$.
For ${\beta} \in {\mathbb N}_0^{d+1}$ with $|{\beta}| \le N$,
\begin{gather}
\label{eq:Q-genfunc}
R_{{\beta},N}^{\kappa} (y): = |y|^N R_{\beta}^{\kappa} \left(\frac{y'}{|y|} \right)
=\frac{(-1)^{|{\beta}|}}{(-N)_{|{\beta}|}}
\sum\limits_{|{\alpha}|=N}
\frac{ N!}{{\alpha}!}{\mathsf Q}_{\beta}({\alpha}; {\kappa}, N) y^{\alpha}.
\end{gather}
\end{thm}

\begin{proof}
Recall that $P_{\nu,N}^{\kappa}(y) = |y|^N P_\nu(y/|y|)$.
Let $n = |{\beta}|$.
It follows from~\eqref{eq:Ra=Pnu},~\eqref{Hahngenfunc} and~\eqref{eq:xi=phi} that
\begin{gather*}
R_{{\beta},N}^{\kappa}(y)=\frac{({\kappa}+{\mathbf{1}})_{\beta}}{({\lambda}_{\kappa})_n}
\sum\limits_{|\nu| =n}\frac{{\mathsf H}_\nu({\beta}; {\kappa},n)}{{\mathsf B}_\nu({\kappa},N)}P_{\nu,N}(y)
\\
\phantom{R_{{\beta},N}^{\kappa}(y)}
=\frac{({\kappa}+{\mathbf{1}})_{\beta}}{ ({\lambda}_{\kappa})_n}
\sum\limits_{|{\alpha}| =N}\frac{N!}{{\alpha}!}
\sum\limits_{|\nu| =n}\frac{{\mathsf H}_\nu({\beta}; {\kappa},n)}{{\mathsf B}_\nu({\kappa},N)}{\mathsf H}_\nu({\alpha};{\kappa},N) y^{\alpha}
\\
\phantom{R_{{\beta},N}^{\kappa}(y)}
= \frac{({\kappa}+{\mathbf{1}})_{\beta}}{ ({\lambda}_{\kappa})_n}
\sum\limits_{|{\alpha}| =N}\frac{N!}{{\alpha}!}
\frac{(-1)^n({\lambda}_{\kappa})_n}{(-N)_n ({\kappa}+{\mathbf{1}})_{\beta}}{\mathsf Q}_{\beta}({\alpha};{\kappa},N)y^{\alpha},
\end{gather*}
which simplif\/ies to the~\eqref{eq:Q-genfunc}.
\end{proof}

{\sloppy The generating function relation~\eqref{eq:Q-genfunc} can be used to derive an explicit expansion of ${\mathsf Q}_{\beta}(\cdot; {\kappa}, N)$.
The following simple lemma is useful.

}

\begin{lem}
For $y \in {\mathbb Z}_N^{d+1}$ and ${\alpha} \in {\mathbb N}_0^{d+1}$ with $|{\alpha}| \le N$,
\begin{gather}
\label{eq:monic-monic}
\frac1{ (-N)_{|{\alpha}|} }
\sum\limits_{|{\gamma}| = N}
\frac{N!}{{\gamma}!}(-{\gamma})_{\alpha} y^{\gamma} =|y|^{N-|{\alpha}|} y^{\alpha}.
\end{gather}
\end{lem}

\begin{proof}
Using the multinomial identity we obtain that
\begin{gather*}
\sum\limits_{|{\gamma}| = N}\frac{N!}{{\gamma}!}(-{\gamma})_{\alpha} y^{\gamma}
= (-1)^{|{\alpha}|}\sum\limits_{|{\gamma}| = N}\frac{N!}{({\gamma}-{\alpha})!} y^{\gamma}
\\
\phantom{\sum\limits_{|{\gamma}| = N}\frac{N!}{{\gamma}!}(-{\gamma})_{\alpha} y^{\gamma}}
= (-1)^{|{\alpha}|}
\sum\limits_{|{\gamma}| = N-|{\alpha}|}
\frac{N!}{{\gamma}!} y^{{\gamma}+{\alpha}} = \frac{(-1)^{|{\gamma}|} N!}{(N-|{\alpha}|)!} |y|^{N-|{\alpha}|} y^{\alpha},
\end{gather*}
which proves the stated identity.
\end{proof}

\begin{thm}
\label{prop:Q=m}
For ${\alpha} \in {\mathbb N}_0^{d+1}$,
\begin{gather}
\label{eq:Q-3F2}
{\mathsf Q}_{\alpha}(x; {\kappa}, N) = \frac{(-N)_{|{\alpha}|}
({\kappa}+{\mathbf{1}})_{\alpha}}{(|{\kappa}|+d+|{\alpha}|)_{|{\alpha}|}}
\sum\limits_{{\gamma} \le {\alpha}}
\frac{(-{\alpha})_{\gamma} (|{\kappa}|+d+|{\alpha}|)_{|{\gamma}|} (-1)^{|{\gamma}|}}{ {\gamma}!
({\kappa}+{\mathbf{1}})_{\gamma} (-N)_{|{\gamma}|}} \mathrm{m}_{\alpha}(x).
\end{gather}
\end{thm}

\begin{proof}
By the expansion of $R_{\alpha}^{\kappa}$ and~\eqref{eq:monic-monic},
\begin{gather*}
|y|^N R_{\alpha}(y')
=\frac{(-1)^{|{\alpha}|} ({\kappa} + {\mathbf{1}})_{\alpha}}{(|{\kappa}|+d+|{\alpha}|)_{|{\alpha}|}}
\sum\limits_{{\gamma} \le {\alpha}}
\frac{(-{\alpha})_{\gamma} (|{\kappa}|+d+|{\alpha}|)_{|{\gamma}|} }{ {\gamma}! ({\kappa}+{\mathbf{1}})_{\gamma}}
|y|^{N-|{\gamma}|} y^{\gamma}
\\
\phantom{|y|^N R_{\alpha}(y')}
=\sum\limits_{|{\beta}|=N}
\frac{N!}{{\beta}!} \left[ \frac{(-1)^{|{\alpha}|} ({\kappa} + {\mathbf{1}})_{\alpha}}{(|{\kappa}|+d+|{\alpha}|)_{|{\alpha}|}}
\sum\limits_{{\gamma} \le {\alpha}}
\frac{(-{\alpha})_{\gamma} (|{\kappa}|+d+|{\alpha}|)_{|{\gamma}|} (-{\beta})_{\gamma}}{{\gamma}!
({\kappa}+{\mathbf{1}})_{\gamma} (-N)_{|{\gamma}|}} \right] X^{\beta}.
\end{gather*}
Comparing with~\eqref{eq:Q-genfunc} proves~\eqref{eq:Q-3F2}.
\end{proof}

\begin{prop}
\label{prop:m=Q}
For $x \in {\mathbb Z}_N^{d+1}$ and ${\alpha} \in {\mathbb N}_0^{d+1}$ with $|{\alpha}| \le N$,
\begin{gather}
\label{eq:monic-Q}
\mathrm{m}_{\alpha}(x)= (-1)^{|{\alpha}|}
\sum\limits_{{\beta} \le {\alpha}}
\frac{(-{\alpha})_{\beta} (-N)_{|{\alpha}|}({\kappa}+{\mathbf{1}})_{\alpha} ({\lambda}_{\kappa})_{2 |{\beta}|}}{
{\beta}! (-N)_{|{\beta}|}({\kappa}+{\mathbf{1}})_{\beta} ({\lambda}_{\kappa})_{|{\alpha}|+ |{\beta}|} }{\mathsf Q}_{\beta}(x; {\kappa}, N),
\end{gather}
where ${\lambda}_{\kappa} = |{\kappa}| +d+1$ as before.
\end{prop}

\begin{proof}
Using the expansion of $Y^{\alpha}$ in terms of $R_{\beta}^{\kappa}$ and the generating function
relation~\eqref{eq:Q-genfunc}, we obtain, for $x \in {\mathbb R}^d$ and $X = (x, 1-|x|)$, that
\begin{gather*}
X^{\alpha}
=({\kappa}+{\mathbf{1}})_{\alpha}\sum\limits_{{\beta} \le {\alpha}}\sum\limits_{{\beta} \le {\alpha}}
\frac{(-1)^{|{\beta}|}(-{\alpha})_{\beta} ({\lambda}_{\kappa})_{2 |{\beta}|}}{{\beta}!
({\kappa}+{\mathbf{1}})_{\beta}({\lambda}_{\kappa})_{|{\alpha}|+|{\beta}|}} R_{\beta}^{\kappa}(y)
\\
\phantom{X^{\alpha}}
=\sum\limits_{|{\gamma}|=N}
\frac{N!}{{\gamma}!} \left[\sum\limits_{{\beta} \le {\alpha}}
\frac{(-{\alpha})_{\beta}({\kappa}+{\mathbf{1}})_{\alpha} ({\lambda}_{\kappa})_{2 |{\beta}|}}
{{\beta}!(-N)_{|{\beta}|}({\kappa}+{\mathbf{1}})_{\beta} ({\lambda}_{\kappa})_{|{\alpha}|+|{\beta}|}}
{\mathsf Q}_{\beta}(x; {\kappa}, N) \right] X^{\gamma}.
\end{gather*}
Setting $X= y /|y|$ with $y \in {\mathbb R}^{d+1}$, we obtain an expansion of $|y|^{N-|{\alpha}|} y^{\alpha}$ in terms
of $y^{\gamma}$, which implies, when comparing with~\eqref{eq:monic-monic}, the identity~\eqref{eq:monic-Q}.
\end{proof}

\begin{defn}
For $n =0,1,\ldots, N$, let $\mathrm{proj}_{{\mathcal V}_n^d({\mathsf H}_{{\kappa},N})}$ denote the orthogonal projection
from $\Pi_N^d$ onto ${\mathcal V}_n^d({\mathsf H}_{{\kappa},N})$.
For ${\alpha} \in {\mathbb N}_0^{d+1}$ with $n \le |{\alpha}| \le N$, define
\begin{gather*}
{\mathsf Q}_{{\alpha},n}
(x; {\kappa}, N): = \frac1{({\kappa}+{\mathbf{1}})_{\alpha}} \mathrm{proj}_{{\mathcal V}_n^d({\mathsf H}_{{\kappa},N})} \mathrm{m}_{\alpha}(x),
\qquad
x \in {\mathbb Z}_N^{d+1}.
\end{gather*}
\end{defn}

If $|{\alpha}| = n$, then ${\mathsf Q}_{{\alpha},n}(\cdot; {\kappa}, N) ={\mathsf Q}_{{\alpha}}
(\cdot; {\kappa}, N)/({\kappa}+{\mathbf{1}})_{\alpha}$.
For $|{\alpha}| \ge n$, it follows from~\eqref{eq:monic-Q} that ${\mathsf Q}_{{\alpha},n}(x; {\kappa}, N)$ has the following expansion:

\begin{prop}
\label{prop:Qan=Q}
For ${\alpha} \in {\mathbb N}_0^{d+1}$ with $|a| \le N$ and $n = 0,1,\ldots, N$,
\begin{gather}
\label{eq:Qan}
{\mathsf Q}_{{\alpha},n}
(x; {\kappa}, N) = \frac{ (-1)^{|{\alpha}|}(-N)_{|{\alpha}|} ({\lambda}_{\kappa})_{2 n}}{
(-N)_{n}({\lambda}_{\kappa})_{|{\alpha}|+n} }
\sum\limits_{|{\beta}|=n}
\frac{(-{\alpha})_{\beta}}{{\beta}! ({\kappa}+{\mathbf{1}})_{\beta}}{\mathsf Q}_{\beta}(x; {\kappa}, N).
\end{gather}
\end{prop}

Evidently, for f\/ixed $n$ and $m \ge n$, the cardinality of the set $\{{\mathsf Q}_{{\alpha},n}
(\cdot; {\kappa}, N): |{\alpha}| =m\}$ is much larger than the dimension of the ${\mathcal V}_n^d({\mathsf H}_{{\kappa},N})$.
Hence, there are ample redundancy in the set.
We will need explicit formulas for ${\mathsf Q}_{{\alpha},n}(x;{\kappa},N)$ for which the following function def\/ined in~\cite{X13} is useful.

\begin{defn}
Let ${\kappa} \in {\mathbb R}^{d+1}$ with ${\kappa}_i > -1$.
For $x, y \in {\mathbb Z}_N^{d+1}$, and $n = 0,1,\ldots$, define ${\mathcal E}_0(x,y; {\kappa})=1$ and
\begin{gather*}
{\mathcal E}_n(x,y; {\kappa}):=\sum\limits_{|{\gamma}| = n}\frac{(-x)_{\gamma} (-y)_{\gamma}}{{\gamma}! ({\kappa}+{\mathbf{1}})_{\gamma}},
\qquad
n =1,2,\dots, N.
\end{gather*}
\end{defn}

\begin{prop}
\label{prop:Q=CE}
For $m \ge n$, ${\alpha} \in {\mathbb Z}_m^{d+1}$ and $x \in {\mathbb Z}_N^{d+1}$,
\begin{gather}
{\mathsf Q}_{{\alpha},n}(x;{\kappa},N)
=\frac{(-N)_m (-m)_n (|{\kappa}|+d+1)_{n} (|{\kappa}|+d+2n)}{n! (|{\kappa}|+d+1)_{m+n}(|{\kappa}|+d+n)}\nonumber
\\
\phantom{{\mathsf Q}_{{\alpha},n}(x;{\kappa},N)=}
\times
\sum\limits_{k=1}
^n \frac{(-n)_k (|{\kappa}|+d+n)_k}{(-m)_k (-N)_k}{\mathcal E}_k({\alpha}, x; {\kappa}).\label{eq:QinE}
\end{gather}
\end{prop}

\begin{proof}
Using ${\lambda}_{\kappa} = |{\kappa}| +d+1$, it follows from~\eqref{eq:Qan} and~\eqref{eq:Q-3F2} that
\begin{gather*}
{\mathsf Q}_{{\alpha},n}
(x;{\kappa},N) = \frac{(-1)^m(-N)_m ({\lambda}_{\kappa})_{2n}}{({\lambda}_{\kappa})_{m+n} (|{\kappa}|+d+n)_n }
\sum\limits_{|{\beta}| = n}
\frac{(-{\alpha})_{\beta}}{{\beta}!}
\sum\limits_{{\gamma} \le {\beta}}
\frac{(-{\beta})_{\gamma} (|{\kappa}|+d+n)_{|{\gamma}|}\mathrm{m}_{\gamma}(x)}{{\gamma}!({\kappa}+{\mathbf{1}})_{\gamma}
(-N)_{|{\gamma}|}}.
\end{gather*}
Since $(-{\gamma})_{\beta} =0$ whenever ${\gamma} > {\beta}$, we can write the sum over ${\gamma} \le {\beta}$ as
${\gamma} \in {\mathbb N}_0^{d+1}$, so that we can consider the sum over $|{\beta}| = n$ f\/irst.
By the multinomial identity,
\begin{gather*}
\sum\limits_{|{\beta}|=n}\frac{(-{\alpha})_{\beta} (-{\beta})_{\gamma}}{{\beta}!}
=\sum\limits_{|{\beta}|=n}\frac{(-{\alpha})_{\beta} (-1)^{|{\gamma}|}}{ ({\beta}-{\gamma})!}
=(-1)^{|{\gamma}|}\sum\limits_{|{\beta}| = n-|{\gamma}|}\frac{(-{\alpha})_{{\beta}+{\gamma}}}{{\beta}!}
\\
\phantom{\sum\limits_{|{\beta}|=n}\frac{(-{\alpha})_{\beta} (-{\beta})_{\gamma}}{{\beta}!}}
= (-1)^{|{\gamma}|} (-{\alpha})_{\gamma}
\sum\limits_{|{\beta}| = n-|{\gamma}|}
\frac{(-{\alpha}+{\gamma})_{{\beta}}}{{\beta}!} = (-1)^{|{\gamma}|} (-{\alpha})_{\gamma}
\frac{(-|{\alpha}|+|{\gamma}|)_{n - |{\gamma}|}}{(n-|{\gamma}|)!}
\\
\phantom{\sum\limits_{|{\beta}|=n}\frac{(-{\alpha})_{\beta} (-{\beta})_{\gamma}}{{\beta}!}}
= \frac{(-{\alpha})_{\gamma} (-|{\alpha}|)_n (-n)_{|{\gamma}|} }{(-|{\alpha}|)_{|{\gamma}|} n!}.
\end{gather*}
Consequently, we obtain
\begin{gather*}
{\mathsf Q}_{{\alpha},n}
(x;{\kappa},N) = \frac{(-1)^m(-N)_m ({\lambda}_{\kappa})_{2n} (-m)_n}{({\lambda}_{\kappa})_{m+n} (|{\kappa}|+d+n)_n
n!}
\sum\limits_{{\gamma}}
\frac{(-n)_{|{\gamma}|} (|{\kappa}|+d+n)_{|{\gamma}|}(-{\alpha})_{\gamma} \mathrm{m}_{\gamma}(x)}{(-m)_{|{\gamma}|}
(-N)_{|{\gamma}|}{\gamma}! ({\kappa}+{\mathbf{1}})_{\gamma} }.
\end{gather*}
Simplifying the constant in front and writing the summation over ${\gamma}$ as $\sum\limits_{k=0}^n\sum\limits_{|{\gamma}| =k}$,
we then obtain~\eqref{eq:QinE}.
\end{proof}

\begin{prop}
\label{prop:Qan=H}
For $m \ge n$, ${\alpha} \in {\mathbb Z}_m^{d+1}$ and $\nu \in {\mathbb N}_0^{d}$ with $|\nu| =n$,
\begin{gather}
\label{Qan=sH}
{\mathsf Q}_{{\alpha},n}
(x; {\kappa}, N) = \frac{(-1)^m(-m)_n (-N)_m ({\lambda}_{\kappa})_{N+n}}{
(-N)_n({\lambda}_{\kappa})_N({\lambda}_{\kappa})_{m+n}}
\sum\limits_{|\nu|=n}
\frac{
{\mathsf H}_\nu({\alpha}; {\kappa}, m)}{
{\mathsf B}_\nu({\kappa},N)}{\mathsf H}_\nu(x;{\kappa},N).
\end{gather}
Conversely,
\begin{gather}
\label{H=Qan}
{\mathsf H}_\nu(x;{\kappa},N) = \frac{(-N)_m}{m!}
\sum\limits_{|{\alpha}|= m}
\frac{({\kappa}+{\mathbf{1}})_{\alpha}}{{\alpha}!}{\mathsf H}_\nu({\alpha}; {\kappa}, m) {\mathsf Q}_{{\alpha},n}
(x;{\kappa},N).
\end{gather}
\end{prop}

\begin{proof}
From~\eqref{eq:Qan} and~\eqref{eq:xi=phi}, it follows that
\begin{gather*}
{\mathsf Q}_{{\alpha},n}
(x; {\kappa}, N) = \frac{ (-1)^{m+n}(-N)_{m} ({\lambda}_{\kappa})_{2 n}}{ ({\lambda}_{\kappa})_{m+n}
({\lambda}_{\kappa})_n}
\sum\limits_{|{\beta}|=n}
\frac{(-{\alpha})_{\beta} }{{\beta}! }
\sum\limits_{|\nu|=n}
\frac{
{\mathsf H}_\nu({\beta}; {\kappa},n)}{
{\mathsf B}_\nu({\kappa},N)}{\mathsf H}_\nu(x; {\kappa},N).
\end{gather*}
Exchanging the order of summations and applying~\eqref{eq:phi-monic}, we obtain~\eqref{Qan=sH} after simplifying the constants.
Conversely, we use the orthogonality of ${\mathsf H}_\nu(\cdot; {\kappa}, m)$ and deduce from~\eqref{Qan=sH} that
\begin{gather*}
\sum\limits_{|{\alpha}| =m}
\frac{({\kappa}+{\mathbf{1}})_{\alpha}}{{\alpha}!}{\mathsf H}_\nu({\alpha}; {\kappa},m) {\mathsf Q}_{{\alpha},n}
(x;{\kappa},N)
\\
\qquad
= \frac{(-m)_n (-N)_m ({\lambda}_{\kappa})_{N+n}}{ (-N)_n({\lambda}_{\kappa})_N({\lambda}_{\kappa})_{m+n}}
\frac{({\lambda}_k)_m}{m!} \frac{
{\mathsf B}_\nu({\kappa}, m)}{
{\mathsf B}_\nu({\kappa}, N)}{\mathsf H}_{\nu}
(x;{\kappa},N).
\end{gather*}
Simplifying the constant by~\eqref{eq:B-A} proves~\eqref{H=Qan}.
\end{proof}

In particular,~\eqref{H=Qan} shows that the set $\{
{\mathsf Q}_{{\alpha},n}
(\cdot;{\kappa},N): |{\alpha}| = m\}$, which is a~subset of the space ${\mathcal V}_n^d({\mathsf H}_{{\kappa},N}
)$, spans the space.

Let ${\mathsf P}_n({\mathsf H}_{{\kappa},N}
; \cdot, \cdot)$ denote the reproducing kernel of ${\mathcal V}_n^d({\mathsf H}_{{\kappa},N}
)$, which is characterized by the requirement that it is an element of ${\mathcal V}_n^d({\mathsf H}_{{\kappa},N}
)$ in either its variable and
\begin{gather*}
\left \langle
{\mathsf P}_n({\mathsf H}_{{\kappa},N}
; x, \cdot),
{\mathsf P} \right \rangle_{{\mathsf H}_{{\kappa},N}
} =
{\mathsf P}(x),
\qquad
\forall\, {\mathsf P} \in {\mathcal V}_n^d({\mathsf H}_{{\kappa},N}
).
\end{gather*}
In terms of the mutually orthogonal basis $\{{\mathsf H}_\nu(\cdot; {\kappa}, N): |\nu| =n\}$, the reproducing kernel can be written as
\begin{gather}
\label{eq:reprod-sH}
{\mathsf P}_n({\mathsf H}_{{\kappa},N};x,y)
=\sum\limits_{|\nu| = n}\frac{{\mathsf H}_\nu(x; {\kappa}, N){\mathsf H}_\nu(y; {\kappa}, N)}{{\mathsf B}_\nu({\kappa}, N)}.
\end{gather}
Our next result shows that this kernel can be expanded in ${\mathsf Q}_{{\alpha},n}
(\cdot; {\kappa}, N)$.

\begin{thm}
\label{thm:reprod-sQ}
For $m \ge n$, ${\alpha} \in {\mathbb Z}_m^{d+1}$ and $x, y \in {\mathbb Z}_N^{d+1}$,
\begin{gather}
\label{eq:reprod-sQ}
{\mathsf P}_n({\mathsf H}_{{\kappa},N}; x,y)
={\mathsf D}_n(m,N) \sum\limits_{|{\alpha}|=m}
\frac{({\kappa}+{\mathbf{1}})_{\alpha}}{{\alpha}!}{\mathsf Q}_{{\alpha},n}(x; {\kappa}, N){\mathsf Q}_{{\alpha},n}(y; {\kappa}, N),
\end{gather}
where
\begin{gather*}
{\mathsf D}_n(m,N): = \frac{(-N)_n m! (|{\kappa}|+d+1)_N (|{\kappa}|+d+1)_{m+n}
}{(-m)_n [(-N)_m]^2 (|{\kappa}|+d+1)_{N+n}}.
\end{gather*}
\end{thm}

\begin{proof}
Let $C_n(m,N)$ denote the constant in from of the sum in~\eqref{Qan=sH}.
Then
\begin{gather*}
\sum\limits_{|{\alpha}| = m}
\frac{({\kappa}+{\mathbf{1}})_{\alpha}}{{\alpha}!}{\mathsf Q}_{{\alpha},n}(x; {\kappa}, N){\mathsf Q}_{{\alpha},n}(y; {\kappa}, N)
\\
\qquad
=[C_n(m,N)]^2\sum\limits_{|\nu|=n}\sum\limits_{|\mu|=n}
\frac{{\mathsf H}_\nu(x;{\kappa},N){\mathsf H}_\mu(y;{\kappa},N)}{{\mathsf B}_\nu({\kappa},N){\mathsf B}_\mu({\kappa},N)}
\!
\sum\limits_{|{\alpha}| = m}
\frac{({\kappa}+{\mathbf{1}})_{\alpha}}{{\alpha}!}{\mathsf H}_\mu({\alpha}; {\kappa}, m){\mathsf H}_\nu({\alpha}; {\kappa}, m)
\\
\qquad
= C_n(m,N)^2 \sum\limits_{|\nu|=n}
\frac{{\mathsf H}_\nu(x;{\kappa},N){\mathsf H}_\nu(y;{\kappa},N)}{{\mathsf B}_\nu({\kappa},N)}
\frac{ ({\lambda}_k)_m{\mathsf B}_\nu({\kappa},m)}{m!{\mathsf B}_\nu({\kappa},N)}
\\
\qquad
= C_n(m,N)^2 \frac{ ({\lambda}_k)_m ({\lambda}_{\kappa})_{m+n}({\lambda}_{\kappa})_N (-N)_n}{m!
({\lambda}_{\kappa})_{N+n}({\lambda}_{\kappa})_m (-m)_n}{\mathsf P}_n({\mathsf H}_{{\kappa},N}; x,y),
\end{gather*}
where we have used~\eqref{eq:B-A}.
This proves~\eqref{eq:reprod-sQ} after simplifying the constants.
\end{proof}

While~\eqref{eq:reprod-sH} follows directly from the def\/inition,~\eqref{eq:reprod-sQ} is by no means trivial since the
set $\{
{\mathsf Q}_{{\alpha}, n}
(\cdot; {\kappa}, N): {\alpha} \in {\mathbb Z}_N^{d+1}\}$ is linearly dependent and is heavily redundant.
In fact, the identity~\eqref{eq:reprod-sQ} is the key result for our discussion on the tight frames in Section~\ref{Sec5}.

\section{Monic Krawtchouk polynomials of several variables}\label{Sec4}
\setcounter{equation}{0}

For each ${\alpha} \in {\mathbb Z}_N^{d+1}$, we def\/ine the monic Krawtchouk polynomial as the orthogonal polynomial that
has $\mathrm{m}_{\alpha}$ as its leading term.

\begin{defn}
For ${\alpha} \in {\mathbb N}_0^{d+1}$ and $|{\alpha}| \le N$, the monic Krawtchouk polynomial
${\mathsf L}_{\alpha}(\cdot; {\kappa}, N) \in {\mathcal V}_{|{\alpha}|}^d({\mathsf K}_{\rho,N})$ is defined uniquely by
\begin{gather*}
{\mathsf K}_{\alpha}(x; \rho, N):= \mathrm{m}
_{\alpha}(x) + q_{\alpha}(x),
\qquad
q \in \Pi_{|{\alpha}|-1}^d,
\qquad
x \in {\mathbb Z}_N^{d+1}.
\end{gather*}
\end{defn}

The properties of such polynomials can be derived from the monic Hahn polynomials\linebreak %
${\mathsf Q}_\nu(\cdot; {\kappa}, N)$.
Indeed, it is easy to see that
\begin{gather*}
\lim_{t \to \infty}t^{-N}{\mathsf H}_{\nu, t {\boldsymbol{\rho}}}(x)={\mathsf K}_{\rho}(x),
\qquad
x \in {\mathbb Z}_N^{d+1},
\qquad
\text{and}
\qquad
\lim_{t \to \infty}{\mathsf B}_{\nu}
(t {\boldsymbol{\rho}}, N) =
{\mathsf C}_{\nu}
(\rho, N).
\end{gather*}
From the f\/irst identity follows readily that the inner product ${\langle} f, g{\rangle}_{
{\mathsf H}_{{\kappa},N}
}$ becomes ${\langle} f, {\gamma}{\rangle}_{
{\mathsf K}_{\rho,N}
}$ if we set ${\kappa} = t {\boldsymbol{\rho}}$ and let $t\to \infty$.
Consequently, we conclude that

\begin{prop}
For $\rho \in (0,1)^d$ with $|\rho| < 1$, $\nu \in {\mathbb N}_0^d$ and $x \in {\mathbb Z}_N^{d+1}$,
\begin{gather*}
\lim_{t\to \infty}{\mathsf Q}_\nu (x; t {\boldsymbol{\rho}}, N) = {\mathsf L}_\nu(x; \rho, N).
\end{gather*}
\end{prop}

Together with the limit relation~\eqref{Hahn-Kraw}, we deduce following relations from Proposition~\ref{prop:H=Q=H}.

\begin{prop}
For $x\in {\mathbb Z}_N^{d+1}$, ${\alpha} \in {\mathbb Z}_n^{d+1}$ and $\nu \in {\mathbb N}_0^d$ with $|\nu| =n$,
\begin{gather}
\label{eq:K=L}
{\mathsf K}_\nu(x;\rho,N) = \frac{(-1)^{n}n!}{(-N)_n}\sum\limits_{|{\alpha}| =n}
\frac{{\mathsf K}_\nu({\alpha};\rho,n)}{{\alpha}!}{\mathsf L}_{\alpha}(x;\rho,N),
\end{gather}
and, conversely,
\begin{gather*}
{\mathsf L}_{\alpha}(x;\rho,N) = (-1)^n (-N)_n {\boldsymbol{\rho}}^{\alpha} \sum\limits_{|\nu| =n}
\frac{{\mathsf K}_\nu({\alpha};\rho,n)}{{\mathsf C}_\nu(\rho, n)}{\mathsf K}_\nu(x;\rho,N).
\end{gather*}
\end{prop}

Furthermore, from Theorem~\ref{prop:Q=m} and Proposition~\ref{prop:m=Q}, we deduce the expansion of monic Krawtchouk
polynomials and its converse.

\begin{thm}
For ${\alpha} \in {\mathbb N}_0^{d+1}$,
\begin{gather*}
{\mathsf L}_{\alpha}(x; \rho, N) = (-N)_{|{\alpha}|}{\boldsymbol{\rho}}^{\alpha}\sum\limits_{{\gamma} \le {\alpha}}
\frac{(-{\alpha})_{\gamma} (-1)^{|{\gamma}|}}{ {\gamma}! (-N)_{|{\gamma}|}{\boldsymbol{\rho}}^{\gamma} }\mathrm{m}_{\alpha}(x).
\end{gather*}
Conversely,
\begin{gather*}
\mathrm{m}_{\alpha}(x)= (-1)^{|{\alpha}|} (-N)_{|{\alpha}|}{\boldsymbol{\rho}}^{\alpha}
\sum\limits_{{\beta} \le {\alpha}}
\frac{(-{\alpha})_{\beta}}{ {\beta}! (-N)_{|{\beta}|}{\boldsymbol{\rho}}^{\beta}}{\mathsf L}_{\beta}(x; {\kappa}, N).
\end{gather*}
\end{thm}

We now def\/ine analogues of ${\mathsf Q}_{{\alpha},n}
(x; {\kappa}, N)$ polynomials.

\begin{defn}
For $n =0,1,\ldots, N$, let $\mathrm{proj}_{{\mathcal V}_n^d({\mathsf K}_{\rho,N}
)}$ denote the orthogonal projection from $\Pi_N^d$ onto ${\mathcal V}_n^d({\mathsf K}_{\rho,N})$.
For ${\alpha} \in {\mathbb N}_0^{d+1}$ with $n \le |{\alpha}| \le N$, define
\begin{gather*}
{\mathsf L}_{{\alpha},n}
(x; {\kappa}, N): = \frac1{({\kappa}+{\mathbf{1}})_{\alpha}} \mathrm{proj}_{{\mathcal V}_n^d({\mathsf K}_{{\kappa},N}
)} \mathrm{m}_{\alpha}(x),
\qquad
x \in {\mathbb Z}_N^{d+1}.
\end{gather*}
\end{defn}

If $|{\alpha}| = n$, then ${\mathsf L}_{{\alpha},n}
(\cdot; \rho, N) =
{\mathsf L}_{{\alpha}}
(\cdot; \rho, N)/{\boldsymbol{\rho}}^{\alpha}$.
It follows that
\begin{gather*}
\lim_{t \to \infty} t^{|{\alpha}|}{\mathsf Q}_{{\alpha},n}
(x; {\kappa}, N) =
{\mathsf L}_{{\alpha},n}
(\cdot; \rho, N).
\end{gather*}
Consequently, from Proposition~\ref{prop:Qan=Q}, we deduce the following relation.

\begin{prop}
For ${\alpha} \in {\mathbb N}_0^{d+1}$ with $|a| \le N$ and $n = 0,1,\ldots, N$,
\begin{gather*}
{\mathsf L}_{{\alpha},n}(x; \rho, N) = \frac{ (-1)^{|{\alpha}|}(-N)_{|{\alpha}|}}{ (-N)_{n}}
\sum\limits_{|{\beta}|=n}
\frac{(-{\alpha})_{\beta}}{{\beta}! {\boldsymbol{\rho}}^{\beta}}{\mathsf L}_{\beta}(x; {\kappa}, N).
\end{gather*}
\end{prop}

Furthermore, we derive from Proposition~\ref{prop:Qan=H} the expansion of ${\mathsf L}_\nu(\cdot; \rho, N)$
in terms of mutually orthogonal ${\mathsf K}_\nu(\cdot; \rho, N)$ and its converse.

\begin{prop}
For $m \ge n$, ${\alpha} \in {\mathbb Z}_m^{d+1}$ and $\nu \in {\mathbb N}_0^{d}$ with $|\nu| =n$,
\begin{gather*}
{\mathsf L}_{{\alpha},n}(x; \rho, N) = \frac{(-m)_n (-N)_m }{ (-N)_n}
\sum\limits_{|\nu|=n}\frac{{\mathsf K}_\nu({\alpha}; \rho, m)}{{\mathsf C}_\nu(\rho,N)}{\mathsf K}_\nu(x;\rho,N).
\end{gather*}
Conversely,
\begin{gather*}
{\mathsf K}_\nu(x;\rho,N) = \frac{(-N)_m}{m!}\sum\limits_{|{\alpha}|= m}
\frac{{\boldsymbol{\rho}}^{\alpha}}{{\alpha}!}{\mathsf K}_\nu({\alpha}; \rho, m) {\mathsf L}_{{\alpha},n}(x;\rho,N).
\end{gather*}
\end{prop}

The analogue of the function ${\mathcal E}_n(x,y; {\kappa})$ is def\/ined as follows:

\begin{defn}
Let $\rho \in (0,1)^d$ with $|\rho| < 1$.
For $x, y \in {\mathbb Z}_N^{d+1}$, and $n = 0,1,\ldots$, define ${\mathcal F}_0(x,y; \rho)=1$ and
\begin{gather*}
{\mathcal F}_n(x,y; \rho): =
\sum\limits_{|{\gamma}| = n}
\frac{(-x)_{\gamma} (-y)_{\gamma}}{{\gamma}! {\boldsymbol{\rho}}^{\gamma}},
\qquad
n =1,2,\dots, N.
\end{gather*}
\end{defn}

This function is def\/ined in~\cite{X13}.
It is easy to see that
\begin{gather*}
\lim_{t \to \infty} t^k {\mathcal E}_k (x,y; t {\boldsymbol{\rho}}, N) = {\mathcal F}_k(x,y; \rho, N),
\qquad
x, y \in {\mathbb Z}_N^{d+1}.
\end{gather*}
Consequently, we can deduce from Proposition~\ref{prop:Q=CE} an expansion of ${\mathsf L}(\nu; \rho, N)$.

\begin{prop}
For $m \ge n$, ${\alpha} \in {\mathbb Z}_m^{d+1}$ and $x \in {\mathbb Z}_N^{d+1}$,
\begin{gather*}
{\mathsf L}_{{\alpha},n}(x;\rho,N) = \frac{(-N)_m (-m)_n}{n!}
\sum\limits_{k=1}^n \frac{(-n)_k}{(-m)_k (-N)_k}{\mathcal F}_k({\alpha}, x; \rho).
\end{gather*}
\end{prop}

Finally, let ${\mathsf P}_n({\mathsf K}_{\rho,N}; \cdot,\cdot)$
denote the reproducing kernel of ${\mathcal V}_n^d({\mathsf K}_{\rho,N})$.
In terms of the mutually orthogonal basis ${\mathsf K}_\nu(\cdot;\rho, N)$,
we can write ${\mathsf P}_n({\mathsf K}_{\rho,N}; \cdot,\cdot)$ as
\begin{gather*}
{\mathsf P}_n(K_{\rho,N}; x,y) =\sum\limits_{|\nu| =n}\frac{{\mathsf K}_\nu(x;\rho, N){\mathsf K}_\nu(y;\rho, N)}{{\mathsf C}(\rho,N)}.
\end{gather*}
It can also be expressed in ${\mathsf L}_{{\alpha},n}(\cdot; \rho, N)$ as the following analogue of Theorem~\ref{thm:reprod-sQ} shows.

\begin{thm}
For $m \ge n$, ${\alpha} \in {\mathbb Z}_m^{d+1}$ and $x, y \in {\mathbb Z}_N^{d+1}$,
\begin{gather*}
{\mathsf P}_n({\mathsf K}_{\rho,N}; x,y) = \frac{(-N)_n m! }{(-m)_n [(-N)_m]^2 }\sum\limits_{|{\alpha}| =m}
\frac{{\boldsymbol{\rho}}^{\alpha}}{{\alpha}!}{\mathsf L}_{{\alpha},n}(x; \rho, N){\mathsf L}_{{\alpha},n}(y; \rho, N).
\end{gather*}
\end{thm}

\section{Finite Tight frames}\label{Sec5}
\setcounter{equation}{0}

For $n \le m \le N$, let us consider the set
\begin{gather*}
\Xi_{m, n}({\kappa},N): = \big\{{\mathsf Q}_{{\alpha},n}(\cdot; {\kappa}, N): {\alpha} \in {\mathbb Z}_m^{d+1}\big\}.
\end{gather*}
Clearly $\Xi_{m,n}({\kappa},N)$ is a~subset of ${\mathcal V}_n^d({\mathsf H}_{{\kappa},N})$.
Moreover, by~\eqref{H=Qan}, $\mathrm{span}\, \Xi_{m,n}({\kappa},N) = {\mathcal V}_n^d({\mathsf H}_{{\kappa},N})$.
The cardinality of $\Xi_{m,n}({\kappa},N)$ is much large than the dimension of ${\mathcal V}_n^d({\mathsf H}_{{\kappa},N})$.
In fact,
\begin{gather*}
\frac{\# \Xi_{m,n}({\kappa},N)}{ \operatorname{dim}{\mathcal V}_n^d({\mathsf H}_{{\kappa},N})}
=\frac{\binom{m+d}{m} }{\binom{n+d-1}{n}} = \frac{m}{d} \Big(\frac{m}{n}\Big)^{d-1} \big(1+ {\mathcal O}\big(m^{-1}\big)\big).
\end{gather*}
The following theorem shows that $\Xi_{m,n}({\kappa},N)$ is a~tight frame for ${\mathcal V}_n^d({\mathsf H}_{{\kappa},N})$.

\begin{thm}
\label{thm:frame-H}
Let $m \ge n$.
Then for all
\begin{gather*}
f(x) =
\sum\limits_{|{\alpha}|=m}
\frac{ ({\kappa} +{\mathbf{1}})_{\alpha} }{\alpha!} \big \langle f,
{\mathsf Q}_{{\alpha},n}
(\cdot; {\kappa}, N) \big \rangle_{
{\mathsf H}_{{\kappa}, N}
}{\mathsf Q}_{x,n}
({\alpha}; {\kappa}, N).
\end{gather*}
Furthermore, $\{
{\mathsf Q}_{{\alpha},n}
(\cdot; {\kappa}, N): |{\alpha}| = m\}$ is a~tight frame of ${\mathcal V}_n^d({\mathsf H}_{n,N}
)$,
\begin{gather}
\label{frame-H}
{\langle} f, f{\rangle}_{{\mathsf H}_{{\kappa},N}}
={\mathsf D}_n(m,N) \sum\limits_{|{\alpha}| = m}\frac{({\kappa} + {\mathbf{1}})_{\alpha}}{{\alpha}!} \big \langle f,{\mathsf Q}_{{\alpha},n}
(\cdot; {\kappa}, N) \big \rangle_{{\mathsf H}_{{\kappa}, N}}^2.
\end{gather}
\end{thm}

\begin{proof}
If $f \in {\mathcal V}_n^d({\mathsf H}_{n,N})$,
then $f(x) = {\langle} f,{\mathsf P}_n({\mathsf H}_{{\kappa},N}; \cdot, x) {\rangle}_{{\mathsf H}_{{\kappa},N}}$,
so that the f\/irst identity follows immediately from~\eqref{eq:reprod-sQ}.
Taking the inner product of the f\/irst identity with $f$ gives the second identity.
\end{proof}

A f\/inite tight frame is equivalent to a~tight frame in an Euclidean space.
From Theorem~\ref{thm:frame-H} and connecting coef\/f\/icients of our Hahn polynomials, we can derive Euclidean tight frames
that can be given explicitly.
Recall that ${\mathsf r} (d,n)= \binom{n+d-1}{n}$.
Fixing a~linear order in the set $\{\nu \in {\mathbb N}_0^d: |\nu| =n\}$,
say the lexicographical order, we denote $x\in {\mathbb R}^{{\mathsf r}(d,n)}$ as $x = (x_\nu: |\nu| = n, \nu \in {\mathbb N}_0^d)$.
The usual inner product in ${\mathbb R}^{{\mathsf r}(d,n)}$ then takes the form
\begin{gather*}
(x,y):=\sum\limits_{|\nu| = n}x_\nu y_\nu,
\qquad
x, y \in {\mathbb R}^{{\mathsf r}(d,n)}.
\end{gather*}

\begin{thm}
\label{thm:Rn-H}
Let $m, n = 1,2,\ldots$ with $m \ge n$.
For ${\alpha} \in {\mathbb Z}_m^{d+1}$, define vectors
\begin{gather*}
{\mathsf h}_{{\alpha},n}
=\sqrt{\frac{m! ({\kappa}+{\mathbf{1}})_{\alpha}}{(|{\kappa}|+d+1)_m {\alpha}!}} \left (\frac{{\mathsf H}_\nu({\alpha},{\kappa}, m)}
{\sqrt{{\mathsf B}_\nu({\kappa}, m)}}: |\nu| = n \right) \in {\mathbb R}^{{\mathsf r}(d,n)}.            
\end{gather*}
Then the set ${\mathsf H}(d, n,m, {\kappa}):= \{{\mathsf h}_{{\alpha},n}
: |{\alpha}| =m\}$ is a~tight frame in ${\mathbb R}^{{\mathsf r}(d,n)}$, that is,
\begin{gather*}
\|x\|^2 = (x,x) =
\sum\limits_{|{\alpha}| = m}
(x,{\mathsf h}_{{\alpha},n})^2,
\qquad
\forall\, x\in {\mathbb R}^{{\mathsf r}(d,n)},
\end{gather*}
and the frame has $\#
{\mathsf H}(d,n,m,{\kappa}) = \binom{m+d}{d}$ elements.
\end{thm}

\begin{proof}
For $y \in {\mathbb R}^{{\mathsf r}(d,n)}$ and $x \in {\mathbb Z}_N^{d+1}$, def\/ine
\begin{gather*}
f_y(x) =
\sum\limits_{|\nu| = n}
y_\nu \frac{
{\mathsf H}_\nu(x; {\kappa}, N)}{
\sqrt{{\mathsf B}_\nu({\kappa},N)}
}.
\end{gather*}
By the orthogonality of ${\mathsf H}_\nu(x;{\kappa},N)$, we immediately have ${\langle} f_y, f_y{\rangle}_{{\mathsf H}_{{\kappa},N}
} = (y,y)$.
On the other hand, by~\eqref{Qan=sH}, we obtain
\begin{gather*}
{\langle} f_y,
{\mathsf Q}_{{\alpha},n}
(\cdot; {\kappa},N) {\rangle}_{
{\mathsf H}_{{\kappa},N}
} = C_n(m,N)
\sum\limits_{|\nu| = n}
y_\nu \frac{
{\mathsf H}_\nu({\alpha}; {\kappa}, m)}{
\sqrt{{\mathsf B}_\nu({\kappa},m)}
} \frac {
\sqrt{{\mathsf B}_\nu({\kappa},m)}
}{
\sqrt{{\mathsf B}_\nu({\kappa},N)}
},
\end{gather*}
where $C_n(m,N)$ is the constant in front of the summation in~\eqref{Qan=sH}.
Using~\eqref{eq:B-A}, it follows readily that
\begin{gather*}
{\langle} f_y,
{\mathsf Q}_{{\alpha},n}
(\cdot; {\kappa},N) {\rangle}_{{\mathsf H}_{{\kappa},N}}^2
= C_n(m,N)^2 \frac{({\lambda}_{\kappa})_{m+n}({\lambda}_{\kappa})_{N}(-N)_n}{({\lambda}_{\kappa})_{N+n}({\lambda}_{\kappa})_{m}(-m)_n}
\frac{(|{\kappa}|+d+1)_m {\alpha}!}{m!({\kappa}+{\mathbf{1}})_{\alpha}}(y,{\mathsf h}_{{\alpha},n})^2.
\end{gather*}
Substituting into~\eqref{frame-H} and simplifying the constant, we see that the right hand side of~\eqref{frame-H}
becomes exactly $\sum\limits_{|{\alpha}| =m}
(y,
{\mathsf h}_{{\alpha},n}
)^2$.
\end{proof}

Recall that ${\mathsf P}_n({\mathsf H}_{{\kappa},m}; \cdot,\cdot)$
denotes the reproducing kernel of ${\mathcal V}_n^d({\mathsf H}_{{\kappa},m})$.

\begin{cor}
Let ${\mathsf h}_{{\alpha},n}$ be the vectors defined in the previous theorem.
Then
\begin{gather}
\label{frame-norm}
({\mathsf h}_{{\alpha},n},{\mathsf h}_{{\beta},n})={\mathsf P}_n({\mathsf H}_{{\kappa},m}; {\alpha},{\beta}),
\qquad
{\alpha}, {\beta} \in {\mathbb Z}_m^{d+1}.
\end{gather}
In particular, in terms of the unit vectors $\widetilde{\mathsf h}_{{\alpha},n}={\mathsf h}_{{\alpha},n}/\|{\mathsf h}_{{\alpha},n}\|$,
\begin{gather*}
\|x\|^2 =\sum\limits_{|{\alpha}| = n}{\mathsf P}_n({\mathsf H}_{{\kappa},m}; {\alpha},{\alpha}) (x, \widetilde {\mathsf h}_{{\alpha},n})^2,
\qquad
\forall\, x \in {\mathbb R}^{{\mathsf r}(d,m)}.
\end{gather*}
\end{cor}

\begin{proof}
The identity~\eqref{frame-norm} follows directly from the def\/inition of ${\mathsf h}_{{\alpha},n}$ and~\eqref{eq:reprod-sH}.
\end{proof}

The tight frames in ${\mathsf H}(d,n,m,{\kappa})$ are given by explicit formulas.
Using the expressions ${\mathsf H}_\nu(x;{\kappa},m)$ in~\eqref{eq:Hn-prod}
and ${\mathsf B}_\nu({\kappa},m)$ in~\eqref{eq:Bnu}, the vectors in the frame can be easily computed.
In the case of $d =2$, ${\mathsf H}(2,n,m,{\kappa})$ is a~tight frame of ${\mathbb R}^{n+1}$ with $(m+2)(m+1)/2$ elements.
Below are several examples for $d =2$ and ${\kappa}=0$, in which the vectors in ${\mathsf H}(2,n,m,{\kappa})$ are column of a~matrix,
which we again call ${\mathsf H}(2,n,m,{\kappa})$.

\begin{exam}
For $d =2$ and ${\kappa} =0$, set ${\mathsf H}(n,m) = {\mathsf H}(2,n,m,0)$.
Then
\begin{gather*}
{\mathsf H}(2,2) =\!\left[
\begin{matrix}
\frac1{\sqrt{6}}&-\sqrt{\frac23}& 0&\frac1{\sqrt{6}}&0&0
\\
\frac1{\sqrt{10}}&0&-\sqrt{\frac25}& -\frac1{\sqrt{10}} &\sqrt{\frac25}& 0
\\
\frac1{\sqrt{30}}&\frac1{\sqrt{30}}&-\sqrt{\frac3{10}}&\frac1{\sqrt{30}}&-\sqrt{\frac3{10}} &\sqrt{\frac3{10}}
\end{matrix}
\right]
\end{gather*}
is a~tight frame in ${\mathbb R}^3$ with $6$ elements,
\begin{gather*}
{\mathsf H}(2,3) = \left[
\begin{matrix}
\sqrt{\frac{3}{14}}\!&-\sqrt{\frac{3}{14}}\!&\frac3{\sqrt{42}}\!&
-\sqrt{\frac{3}{14}}\!&-\sqrt{\frac{2}{21}}\!&0&\sqrt{\frac{3}{14}}\!&\frac3{\sqrt{42}}\!&0&0
\\[1.99mm]
\frac3{\sqrt{70}}\!&\frac1{\sqrt{70}}\!&-\frac3{\sqrt{70}}\!&
-\frac1{\sqrt{70}}\!&0&-\frac4{\sqrt{70}}\!&-\frac3{\sqrt{70}}\!&\frac 3{\sqrt{70}}\!&\frac4{\sqrt{70}}\!&0
\\[1.99mm]
\sqrt{\frac{3}{70}}\!&\sqrt{\frac{3}{70}}\!&-\sqrt{\frac{5}{42}}\!&
\sqrt{\frac{3}{70}}\!& -\sqrt{\frac{5}{42}}\!&-\sqrt{\frac{3}{70}}\!&\sqrt{\frac{3}{70}}\!&-\sqrt{\frac{5}{42}}\!&-\sqrt{\frac{3}{70}}\!&\sqrt{\frac{27}{70}}
\end{matrix}
\right]
\end{gather*}
is a~tight frame in ${\mathbb R}^3$ with $10$ elements, and    
\begin{gather*}
{\mathsf H}(3,3) =\!\left[
\begin{matrix}
\frac{-1}{2\sqrt{5}}&\frac3{2\sqrt{5}}&0&-\frac3{2\sqrt{5}}&0&0&\frac1{2\sqrt{5}}&0&0 &0
\\[1.99mm]
\frac{-1}{2\sqrt{7}}&\frac{-1}{2\sqrt{7}}&\frac{-1}{\sqrt{7}}&\frac1{2\sqrt{7}}&\frac{-2}{\sqrt{7}}&0&\frac{-1}{2\sqrt{7}}&\frac1{\sqrt{7}}&0&0
\\[1.99mm]
\frac{-\sqrt{3}}{2\sqrt{35}}&\frac{-1}{2\sqrt{105}}&\sqrt{\frac5{21}}&
\frac1{2\sqrt{105}}&0&-\sqrt{\frac5{21}}&\frac{\sqrt{3}}{2\sqrt{35}}&-\sqrt{\frac5{21}}&\sqrt{\frac5{21}}&0
\\[1.99mm]
\frac{-1}{2\sqrt{35}}&\frac{-1}{2\sqrt{35}}&\frac{2}{\sqrt{35}}&
\frac{-1}{2\sqrt{35}}&\frac2{\sqrt{35}}&-\frac{3}{\sqrt{35}}&\frac{-1}{2\sqrt{35}}&\frac{2}{\sqrt{35}}&\frac{-3}{\sqrt{35}}&\frac{-2}{\sqrt{35}}
\end{matrix}
\right]
\end{gather*}
is a~tight frame in ${\mathbb R}^4$ with $10$ elements.
\end{exam}

For f\/ixed $d$, $n$, $m$, we obtain a~family of tight frames in ${\mathbb R}^{{\mathsf r}(d,n)}$
with $d+1$ parameters ${\kappa}= ({\kappa}_1, \ldots, {\kappa}_{d+1})$.
Although the proof is established for ${\kappa}_i > -1$, the analytic continuation shows that we obtain a~tight frame
for ${\kappa} \in {\mathbb R}^{d+1}$ as long as $h_{{\alpha},n}$ is well def\/ined.
A frame is called normalized tight frame if all elements of the frame have the same norm~\cite{BF}.
Our tight frame is not normalized in general.
In the case of $d =2$ and $m = n =2$, we can choose ${\kappa}$ to obtain a~normalized tight frame.

\begin{exam}
In the case of $d =2$, $m = n =2$, choosing ${\kappa}_1={\kappa}_2={\kappa}_3 = \frac12
\sqrt{-7-\sqrt{17}
}$, then the row vectors in
\begin{gather*}
{\mathsf H}(2,2,2,{\kappa}) =\!\left[
\begin{matrix}
\frac12\sqrt{-3 + \sqrt{17}}&\frac14\sqrt{-1 + \sqrt{17}}& \frac14\sqrt{21 -5 \sqrt{17}}
\\[1.99mm]
\frac {\sqrt{5 - \sqrt{17}}}{\sqrt{2}}&0&\frac12\sqrt{-4 + \sqrt{17}}
\\[1.99mm]
0&-\frac14\sqrt{-3 + \sqrt{17}}&-\frac14\sqrt{-1 + \sqrt{17}}
\\[1.99mm]
\frac12\sqrt{-3 + \sqrt{17}}&-\frac14\sqrt{-1 + \sqrt{17}}& \frac14\sqrt{21 -5 \sqrt{17}}
\\[1.99mm]
0&\frac14\sqrt{-3 + \sqrt{17}}&-\frac14\sqrt{-1 + \sqrt{17}}
\\[1.99mm]
0&0&\sqrt{\frac12}
\end{matrix}
\right]
\end{gather*}
form a~normalized tight frame in ${\mathbb R}^3$ with $6$ vectors.
\end{exam}

For the Krawtchouk polynomials, the analogue of Theorem~\ref{thm:frame-H} is as follows.

\begin{thm}
\label{thm:frame-K}
Let $m \ge n$.
Then for all
\begin{gather*}
f(x) = \frac{(-N)_m (-m)_n}{n!}
\sum\limits_{|{\alpha}|=m}
\frac{ {\boldsymbol{\rho}}^{\alpha} }{\alpha!} \big \langle f,
{\mathsf K}_{{\alpha},n}
(\cdot; \rho, N) \big \rangle_{
{\mathsf K}_{\rho, N}
}{\mathsf K}_{x,n}
({\alpha}; \rho, N).
\end{gather*}
Furthermore, $\{{\mathsf Q}_{{\alpha},n}(\cdot; {\kappa}, N): |{\alpha}| = m\}$ is a~tight frame of ${\mathcal V}_n^d({\mathsf H}_{n,N})$,
\begin{gather*}
{\langle} f, f{\rangle}_{{\mathsf K}_{\rho,N}} =
\frac{(-N)_m (-m)_n}{n!}\sum\limits_{|{\alpha}| = m}\frac{{\boldsymbol{\rho}}^{\alpha}}{{\alpha}!} \big \langle f,
{\mathsf K}_{{\alpha},n}(\cdot; \rho, N) \big \rangle_{{\mathsf K}_{\rho, N}}^2.
\end{gather*}
\end{thm}

From Theorem~\ref{thm:frame-K} and using the connecting relation in~\eqref{eq:K=L}, we obtain the following analogue of
Theorem~\ref{thm:Rn-H} on tight frames in Euclidean spaces.

\begin{thm}
Let $m, n = 1,2,\ldots$ with $m \ge n$.
For ${\alpha} \in {\mathbb Z}_m^{d+1}$, define vectors
\begin{gather*}
{\mathsf k}_{{\alpha},n}
=\sqrt{\frac{m! {\boldsymbol{\rho}}^{\alpha}}{{\alpha}!}} \left (\frac{
{\mathsf K}_\nu({\alpha},\rho, m)}{
\sqrt{ {\mathsf C}_\nu(\rho, m)}
}: |\nu| = n \right) \in {\mathbb R}^{
{\mathsf r}(d,n)}.
\end{gather*}
Then the set ${\mathsf K}(d, n,m, \rho):= \{{\mathsf k}_{{\alpha},n}: |{\alpha}| =m\}$
is a~tight frame in ${\mathbb R}^{{\mathsf r}(d,n)}$, that is,
\begin{gather*}
(x,x) =
\sum\limits_{|{\alpha}| = m}
(x,
{\mathsf k}_{{\alpha},n}
)^2
\end{gather*}
with $\#
{\mathsf K}(d,n,m,\rho) = \binom{m+d}{d}$ elements.
Furthermore, for ${\alpha}, {\beta} \in {\mathbb Z}_N^{d+1}$,
\begin{gather*}
({\mathsf k}_{{\alpha},n},{\mathsf k}_{{\beta},n}) ={\mathsf P}_n({\mathsf K}_{\rho, m}; {\alpha},{\beta}).
\end{gather*}
\end{thm}

This theorem can also be deduced from Theorem~\ref{thm:Rn-H} by setting ${\kappa} = t {\boldsymbol{\rho}}$ and taking
the limit $t \to \infty$.
The limiting process shows that the tight frames ${\mathsf K}(d,n,m,\rho)$ are not included in ${\mathsf H}(d,n,m,{\kappa})$.
In terms of explicit formulas, the Krawtchouk polynomials are simpler than the Hahn polynomials.

{\sloppy The tight frames ${\mathsf K}(d, n,m, \rho)$ are given by explicit formulas.
Using the expression of ${\mathsf K}_\nu(x; \rho,m)$ in~\eqref{eq:KrawK}
and ${\mathsf C}_\nu(\rho,m)$ in~\eqref{eq:C(rho-N)}, the frame elements can be computed easily.
Below are several examples for $d =2$ and $\rho = (\frac13,\frac13)$, in which the vectors in ${\mathsf K}(2,n,m,{\kappa})$
are column of a~matrix, which we again call ${\mathsf K}(2,n,m,{\kappa})$.

}

\begin{exam}
For $d =2$ and ${\kappa} =0$, set ${\mathsf K}(n,m) = {\mathsf K}(2,n,m,1/3,1/3)$.
Then
\begin{gather*}
{\mathsf K}(2,2) =\!\left[
\begin{matrix}
\frac12&-\frac1{\sqrt{2}}&0&\frac12&0&0
\\[1.99mm]
\frac1{\sqrt{6}}&0&-\frac1{\sqrt{3}}&-\frac1{\sqrt{6}}&\frac1{\sqrt{3}}& 0
\\[1.99mm]
\frac16&\frac1{3\sqrt{2}}&-\frac{\sqrt{2}}3&\frac16&-\frac{\sqrt{2}}3&\frac{2}{3}
\end{matrix}
\right]
\end{gather*}
is a~tight frame in ${\mathbb R}^3$ with $6$ elements,
\begin{gather*}
{\mathsf K}(2,3) =\!\left[
\begin{matrix}
\frac12&-\frac1{2\sqrt{3}}&\frac1{2\sqrt{3}}&
-\frac1{2\sqrt{3}}&-\frac1{\sqrt{6}}&0&\frac1{\sqrt{2}} &\frac1{2\sqrt{3}}&0 &0
\\[1.99mm]
\frac1{\sqrt{6}}&\frac1{3\sqrt{2}}&-\frac1{3\sqrt{2}}&
-\frac1{3\sqrt{2}}&0&-\frac{\sqrt{2}}{3}&-\frac1{\sqrt{6}}&\frac1{3\sqrt{2}}& \frac{\sqrt{2}}3&0
\\[1.99mm]
\frac16&\frac1{2\sqrt{3}}&-\frac1{2\sqrt{3}}&
\frac1{2\sqrt{3}}&-\frac1{\sqrt{6}} &0&\frac16&-\frac1{2\sqrt{3}} &0&\frac{2}{3}
\end{matrix}
\right]
\end{gather*}
is a~tight frame in ${\mathbb R}^3$ with $10$ elements, and
\begin{gather*}
{\mathsf K}(3,3) =\!\left[
\begin{matrix}
-\frac1{2\sqrt{2}}&\frac{\sqrt{3}}{2\sqrt{2}}&0&-\frac{\sqrt{3}}{2\sqrt{2}}&0&0&\frac1{2\sqrt{2}}&0&0 &0
\\[1.99mm]
-\frac1{2\sqrt{2}}&-\frac1{2\sqrt{6}}&-\frac1{\sqrt{6}}&\frac1{2\sqrt{6}}&-\frac1{\sqrt{3}}&0&-\frac1{2\sqrt{2}}&\frac1{\sqrt{6}}&0&0
\\[1.99mm]
-\frac1{2\sqrt{6}}&-\frac1{6\sqrt{2}}&-\frac{\sqrt{2}}3&\frac1{6\sqrt{2}}&0&-\frac{\sqrt{2}}3&\frac1{2\sqrt{6}}&
-\frac{\sqrt{2}}3&-\frac{\sqrt{2}}3&0
\\[1.99mm]
-\frac1{6\sqrt{6}}&-\frac1{6\sqrt{2}}&\frac1{3\sqrt{2}}&-\frac1{6\sqrt{2}}&\frac13&-\frac{\sqrt{2}}3&-\frac1{6\sqrt{6}}&\frac1{3\sqrt{2}}&
-\frac{\sqrt{2}}{3}&\frac{2\sqrt{2}}{3\sqrt{3}}
\end{matrix}
\right]
\end{gather*}
is a~tight frame in ${\mathbb R}^4$ with $10$ elements.
\end{exam}

Using the tight frames of ${\mathcal V}_n^d({\mathsf H}_{{\kappa},N})$ as building blocks,
we can also build tight frame for the space $\Pi_N^d$ of polynomials of degree at most $N$ in
$d$-variables under the inner product ${\langle} \cdot,\cdot{\rangle}_{{\mathsf H}_{{\kappa},N}}$.
For $1 \le j \le N$, let $m_j$ be positive integers such that $j \le m_j \le N$.
Let ${\mathbf m} = \{m_1,\ldots, m_N\}$.
Def\/ine
\begin{gather*}
\Xi_{{\mathbf m}}({\kappa},N):= \{ {\mathbf{1}} \} \bigcup_{n=1}^N \Xi_{m_n, n}({\kappa},N).
\end{gather*}

\begin{thm}
\label{thm:frame-N-H}
For every polynomial $f \in \Pi_N^d$ of degree at most $N$ in $d$-variables,
\begin{gather*}
f = \big \langle f, 1 \big \rangle_{
{\mathsf H}_{{\kappa}, N}
} +
\sum\limits_{n=1}
^N
{\mathsf D}_n(m_n,N) \sum\limits_{|{\alpha}|=m_n}
\frac{ ({\kappa} +{\mathbf{1}})_{\alpha} }{\alpha!} \big \langle f,
{\mathsf Q}_{{\alpha},n}
(\cdot; {\kappa}, N) \big \rangle_{
{\mathsf H}_{{\kappa}, N}
}{\mathsf Q}_{x,n}
({\alpha}; {\kappa}, N), 
\end{gather*}
and, in particular,
\begin{gather}
\label{frame-N-H}
{\langle} f, f{\rangle}_{
{\mathsf H}_{{\kappa},N}
} ={\langle} f, 1{\rangle}_{
{\mathsf H}_{{\kappa},N}
}^2 +
\sum\limits_{n=1}
^N
{\mathsf D}_n(m_n,N) \sum\limits_{|{\alpha}|=m_n}
\frac{({\kappa} + {\mathbf{1}})_{\alpha}}{{\alpha}!} \big \langle f,
{\mathsf Q}_{{\alpha},n}
(\cdot; {\kappa}, N) \big \rangle_{
{\mathsf H}_{{\kappa}, N}
}^2.
\end{gather}
\end{thm}

\begin{proof}
The orthogonal expansion of $f \in \Pi_N^d$ in the Hahn polynomials is given by
\begin{gather*}
f =
\sum\limits_{n=0}
^N \big \langle f,
{\mathsf P}_n({\mathsf H}_{{\kappa},N}
; x, \cdot) \big \rangle_{
{\mathsf H}_{{\kappa},N}
}.
\end{gather*}
Since ${\mathsf P}_0({\mathsf H}_{{\kappa}, N}
; x,y)=1$, the the f\/irst identity follows from from~\eqref{eq:reprod-sQ}, which implies~\eqref{frame-N-H}.
\end{proof}

There is also a~straightforward analogue for the Krawtchouk polynomials, which we shall not state.

Since $\{{\mathsf H}_\nu(\cdot; {\kappa}, N): |\nu| \le N\}$ is a~mutually orthogonal basis, by writing $f$ as
\begin{gather*}
f_y(x) =\sum\limits_{|\nu| \le N}y_\nu{\mathsf H}_\nu(x; {\kappa}, N)
\end{gather*}
and using the connection~\eqref{Qan=sH},  
each tight frame in~\eqref{frame-N-H} corresponds to a~tight frame for ${\mathbb R}^{{\mathsf n}(d,N)}$, where
\begin{gather*}
{\mathsf n}(d, N) = \binom{N+d}{d},
\end{gather*}
and the frame elements can be computed directly in terms of ${\mathsf H}_\nu(x;{\kappa}, m_n)$.
We shall not written down these tight Euclidean frames but turn our attention to another way of constructing tight
Euclidean frames using the Hahn polynomials, this time using ${\mathsf Q}_{{\alpha},n}
(\cdot; {\kappa}, N)$ instead of ${\mathsf H}_\nu(\cdot; {\kappa}, N)$.

This construction is based on the observation that the weight function ${\mathsf H}_{{\kappa},N}$
of the Hahn polynomials becomes $1$ when ${\kappa} = 0$, so that
the inner product ${\langle} \cdot, \cdot {\rangle}_{{\mathsf H}_{0,N}}$ is a~constant multiple of the inner product $(\cdot,\cdot)$
of the Euclidean space ${\mathbb R}^{{\mathsf n}(d,N)}$.
We use Theorem~\ref{thm:frame-N-H} with ${\kappa}=0$ and $m_n =N$ for $n=1,2,\ldots, N$.

\begin{thm}
For ${\alpha} \in {\mathbb Z}_N^{d+1}$, define the vector
\begin{gather}
\label{eq:sq-an}
{\mathsf q}_{{\alpha}, n}
: = \left({\mathsf Q}_{{\alpha},n}
({\beta}; 0, N): {\beta}\in {\mathbb Z}_N^{d+1}\right) \in {\mathbb R}^{
{\mathsf n}(d,N)},
\end{gather}
and define ${\mathsf q}_0: = \sqrt{\tfrac{N!}{(d+1)_N}}{\mathbf{1}}$.
Then the set
\begin{gather*}
\Xi(d,N): = \big\{{\mathsf q}_0 \big\}\cup \big\{{\mathsf q}_{{\alpha},n}: {\alpha} \in {\mathbb Z}_N^{d+1}, n =1,\ldots, N \big\}
\end{gather*}
is a~tight frame of ${\mathbb R}^{{\mathsf n}(d,N)}$ with $\# \Xi(d,N) = 1 + N \times{\mathsf n}(d,N)$ elements, that is,
\begin{gather}
\label{eq:frame-R-Q}
(x,x) = (x,{\mathsf q}_0)^2 + \sum\limits_{n=1}^N\sum\limits_{|{\alpha}| = N}(x,{\mathsf q}_{{\alpha},n})^2.
\end{gather}
\end{thm}

\begin{proof}
For ${\kappa} =0$, the inner product becomes
\begin{gather*}
{\langle} f, g{\rangle}_{
{\mathsf H}_{0,N}
} = \frac{N!}{(d+1)_N}
\sum\limits_{|{\alpha}|=N}
f({\alpha}) g({\alpha}).
\end{gather*}
Since $f \in \Pi_N^d$ is uniquely determined by its values on the points $\{{\alpha}: {\alpha} \in {\mathbb
Z}_N^{d+1}\}$, we see that ${\langle} f, g{\rangle}_{{\mathsf H}_{0,N}} = \frac{N!}{(d+1)_N} ({\mathsf f},{\mathsf g})$,
where ${\mathsf f} = (f({\beta}): {\beta} \in {\mathbb Z}_N^{d+1}) \in {\mathbb R}^{{\mathsf n}(d,N)}$.
For $m = N$ and ${\kappa} =0$, ${\mathsf D}_n(N,N) = (d+1)_N/ N!$ is independent of $n$.
Consequently, we can write~\eqref{frame-N-H}
as~\eqref{eq:frame-R-Q}.
\end{proof}

\begin{cor}
For $n \ge 1$ and ${\mathsf q}_{{\alpha},n}$ in~\eqref{eq:sq-an},
\begin{gather*}
({\mathsf q}_{{\alpha},n},
{\mathsf q}_{{\beta},n}) = \frac{{\mathsf P}_n(H_{0,N}; {\alpha}, {\beta})}{{\mathsf n}(d,N)},
\qquad
{\alpha},{\beta} \in {\mathbb Z}_N^{d+1}.
\end{gather*}
\end{cor}

\begin{proof}
From~\eqref{eq:reprod-sQ} with $m = N$ and ${\kappa} =0$, it follows that
\begin{gather*}
({\mathsf q}_{{\alpha},n},{\mathsf q}_{{\beta},n})
=\sum\limits_{|{\gamma}|=N}{\mathsf Q}_{{\alpha},n}({\gamma}; 0,N){\mathsf Q}_{{\beta},n}({\gamma}; 0,N)
=\frac1{{\mathsf D}_n(N,N)}{\mathsf P}_n(H_{0,N}; {\alpha}, {\beta}).
\end{gather*}
Furthermore, for ${\kappa} =0$, ${\mathsf D}_n(N,N) = (d+1)_N/ N! = \binom{N+d}{N} ={\mathsf n}(d,N)$.
\end{proof}

The elements of the tight frame $\Xi(d,N)$ are given in explicit formulas in terms of ${\mathsf Q}_{{\alpha},n}(\cdot; {\kappa}, N)$.
In the case of $d =2$, we have
\begin{gather*}
{\mathsf Q}_{{\alpha},n}
(x; {\kappa}, N) = \frac{[(-N)_n]^2 (3)_{2n}}{n! (3)_{n+N} (n+2)_n}
\sum\limits_{k=0}^n \frac{(-n)_k (n+2)_k}{[(-N)_k]^2}{\mathcal E}_k ({\alpha},x),
\qquad
{\alpha}, x\in {\mathbb Z}_N^{3},
\end{gather*}
where ${\mathcal E}_0(x,y) = 1$ and, for $k=1,2,\ldots, n$,
\begin{gather*}
{\mathcal E}_k({\alpha},x)
=\sum\limits_{|{\gamma}| = k}\frac{(-{\alpha})_{\gamma} (-x)_{\alpha}}{{\gamma}!^2}
=\sum\limits_{|{\gamma}| = k}\binom{{\alpha}}{{\gamma}} \binom{x}{{\gamma}}.
\end{gather*}
The vectors ${\mathsf q}_{{\alpha},n}$ are given by evaluations of these polynomials.
One interesting feature is that all frame elements other than ${\mathsf q}_0$ have rational entries.
We give the f\/irst two cases as examples.
Let $e_1 = (1,0,0)$, $e_2 = (0,1,0)$, $e_3 = (0,0,1)$.

\begin{exam}
For $d= 2$ and $N =1$, it follows readily that
\begin{gather*}
{\mathsf Q}_{{\alpha},1}
(x,0,1) = \frac13 (1- 3 {\mathcal E}_1({\alpha},x))
\qquad
\text{and}
\qquad
{\mathcal E}_1({\alpha},x) = {\alpha}_1 x_1 +{\alpha}_2 x_2 + {\alpha}_3 x_3.
\end{gather*}
Evaluating at the set $\{{\alpha}: {\alpha} \in {\mathbb Z}_1^3\} = \{e_1, e_2,e_3 \}$ shows that the column of the
matrix
\begin{gather*}
\Xi(2,1)= \left[
\begin{matrix}
\frac1{\sqrt{3}}&\frac23&-\frac13&-\frac13
\\[1.99mm]
\frac1{\sqrt{3}}&-\frac13&\frac23&-\frac13
\\[1.99mm]
\frac1{\sqrt{3}}&-\frac13&-\frac13&\frac23
\end{matrix}
\right]
\end{gather*}
form a~tight frame in ${\mathbb R}^3$ with $4$ vectors.
\end{exam}

\begin{exam}
For $d= 2$ and $N =2$, ${\mathbb Z}_2^3 =\{2e_1,2e_2,2e_3, e_1+e_2,e_2+e_3, e_+e_1\}$.
We have ${\mathsf q}_0 =\frac1{
\sqrt{6}
}(1,1,1,1,1,1)$.
For $n \ge 1$, there are two cases.
For ${\alpha} \in {\mathbb Z}_2^3$ and $n=1$, we obtain
\begin{gather*}
{\mathsf Q}_{{\alpha},1}
(x,0,1) = -\tfrac4{15} \left (1- \tfrac34 {\mathcal E}_1({\alpha},x)\right),
\qquad
{\mathcal E}_1(2e_i,x) = 2x_i,
\qquad
{\mathcal E}_1(e_i+e_j) = x_i+x_j,
\end{gather*}
and for ${\alpha} \in {\mathbb Z}_2^3$ and $n =2$, we obtain
\begin{gather*}
{\mathsf Q}_{{\alpha},2}
(x,0,1) = \tfrac1{10} \left (1- 2 {\mathcal E}_1({\alpha},x) + 10 {\mathcal E}_2({\alpha},x)\right),
\qquad
\end{gather*}
where ${\mathcal E}_1({\alpha},x)$ is as before and
\begin{gather*}
{\mathcal E}_2(2e_i,x) = \tfrac12 x_i(x_i-1),
\qquad
{\mathcal E}_2(2e_i+e_j,x) = x_1x_2.
\end{gather*}
Evaluating these two functions at elements in ${\mathbb Z}_2^3$, we obtain ${\mathsf q}_{{\alpha},1}$ and ${\mathsf q}_{{\alpha},2}$.
The set $\Xi(2,2)$ consists of columns of the matrix
\begin{gather*}
\left[
\begin{array}{ccccccccccccc}
\frac1{\sqrt{6}}&-\frac3{10}&-\frac3{10}&\frac1{10}&\frac1{10}& \frac3{10}& \frac1{10}& -\frac4{15}& \frac{2}{15}&\frac2{15}
& -\frac4{15}&-\frac4{15}&\frac8{15}
\\[1.99mm]
\frac1{\sqrt{6}}&-\frac3{10}&\frac1{10}& -\frac3{10}& \frac1{10}& \frac1{10}& \frac3{10}& \frac2{15}& -\frac{4}{15}&
\frac2{15}&-\frac4{15}& \frac8{15}& -\frac4{15}
\\[1.99mm]
\frac1{\sqrt{6}}&\frac1{10}&-\frac3{10}&-\frac3{10}& \frac3{10}& \frac1{10}& \frac1{10}& \frac2{15}& \frac{2}{15}&
-\frac4{15}&\frac8{15}& -\frac4{15}& -\frac4{15}
\\[1.99mm]
\frac1{\sqrt{6}}&\frac7{10}&-\frac1{10}&-\frac1{10}& \frac1{10}& -\frac3{10}& -\frac3{10}& -\frac1{15}& -\frac1{15}&
\frac2{15}&-\frac4{15}& \frac2{15}& \frac2{15}
\\[1.99mm]
\frac1{\sqrt{6}}&-\frac1{10}&\frac7{10}&-\frac1{10}&-\frac3{10}& -\frac3{10}& \frac1{10}& -\frac1{15}& \frac{2}{15}&
-\frac1{15}&\frac2{15}& -\frac4{15}& \frac2{15}
\\[1.99mm]
\frac1{\sqrt{6}}&-\frac1{10}&-\frac1{10}&\frac7{10}&-\frac3{10}& \frac1{10}& -\frac3{10}& \frac2{15}& -\frac1{15}&
-\frac1{15}&\frac2{15}& \frac2{15}& -\frac4{15}
\end{array}
\right]
\end{gather*}
and these column vectors form a~tight frame in ${\mathbb R}^6$ with 13 vectors.
Apart from the first vector~${\mathsf q}_0$, all other vectors have rational entries.
\end{exam}

For $d =2$ and each $N =1,2,\ldots$, the identity~\eqref{eq:frame-R-Q} gives a~tight frame of ${\mathbb R}^{
{\mathsf n}(2,N)}$ with $M:= N \times
{\mathsf n}(2,N) +1$ vectors.
The next one is a~tight frame of ${\mathbb R}^{10}$ with 31 vectors.
Except ${\mathsf q}_0$, all other frame elements in a~tight frame of this family have rational entries.

\subsection*{Acknowledgements}
The work was supported in part by NSF Grant DMS-1106113.

\pdfbookmark[1]{References}{ref}
\LastPageEnding

\end{document}